\documentclass[11pt]{article}
\usepackage{a4}
\usepackage{graphicx}
\usepackage{parskip}

\usepackage{amsfonts}
\usepackage{natbib}
\usepackage{algorithm}
\usepackage{algorithmic}

\usepackage{latexsym}

\newtheorem{theorem}{Theorem}[section]

\newtheorem{lemma}{Lemma}[section]
\newtheorem{corollary}{Corollary}[section]
\newtheorem{remark}{Remark}[section]

\newtheorem{condition}{Condition}[section]

\newcommand{\R}{\mathbb{R}}

\newcommand{\PP} {{  \rm I\hskip-0.22em P}}

\newcommand{\EE} {{\rm I\hskip-0.48em E}}

\usepackage{a4}

\usepackage{graphicx}

\usepackage{parskip}

\usepackage{amsfonts}

\usepackage{natbib}

\usepackage{algorithm}

\usepackage{algorithmic}

\usepackage{amsmath, amssymb}

\usepackage{latexsym}

\begin{document}

\centerline{\bf A pivotal transform for the high-dimensional location-scale model}

\vskip .1in

\centerline{Sara van de Geer, ETH Z\"urich}

\centerline{Sylvain Sardy, Universit\'e de Gen\`eve}

\centerline{Maxime van Cutsem, Universit\'e de Gen\`eve}

\vskip .1in
\centerline{December 19, 2025}

{\bf Abstract} We study the high-dimensional linear model with noise distribution known up
to a scale parameter. With an $\ell_1$-penalty on the regression coefficients, we show that
a transformation of the log-likelihood allows for  a choice of the tuning parameter
not depending on the scale parameter.  This transformation is a generalization 
of the square root Lasso for quadratic loss.  The tuning parameter can asymptotically be taken
at the detection edge. We establish an oracle inequality, variable selection and 
asymptotic efficiency of the
estimator of the scale parameter and the intercept. The examples include Subbotin distributions and
the Gumbel distribution.

{\it MSC2020 Subject Classification} 62J07, 62J99

{\it Keywords and Phrases} asymptotic efficiency, detection edge, high-dimensional linear model,
location-scale model, oracle inequality, sparsity, variable selection

\section{Introduction}

We study the high-dimensional linear model, with noise distribution known up to a scale parameter.
The density of the noise is assumed to be log-concave, and 
the regression coefficients are assumed to obey a sparsity condition. The variance of the noise then exists, and so one may consider applying the 
square root Lasso (\cite{belloni2011square}, 
\cite{sunzhang11}), based on the least squares loss function, i.e.\ on quadratic loss. We  propose however to use 
a pivotal transformation of minus-log-likelihood loss, which generalizes the square root Lasso to the case of
non-Gaussian noise. We apply the $\ell_1$-penalty on the regression coefficients, with tuning parameter
$\lambda$.

Our aim is threefold.
First of all, we aim at showing that there
is a universal choice of the tuning parameter $\lambda$, which is in particular independent of the unknown scaling
parameter. Secondly, we want that $\lambda$ can be chosen close to the detection edge. The detection edge can be described as follows. Consider the Lasso for the case with known scale parameter and with tuning parameter $\lambda$, given in equation (\ref{Lasso.equation}) below. Suppose the null-model holds, i.e.\  all regression coefficients are zero.
Then,  for $0< \alpha < 1 $, the phase transition at level $1- \alpha$ is the value of the $(1-\alpha)$-quantile $F^{-1} (1- \alpha)$ of a given random variable $\lambda^*$,  such that with 
probability  asymptotically equal to $1- \alpha$,  the Lasso in (\ref{Lasso.equation}) puts all regression
coefficients to zero when $\lambda$ larger than $F^{-1} (1-\alpha) $. The transformed Lasso given 
in (\ref{exp-Lasso.equation}) deals with scale parameter unknown. In 
Section \ref{phase-transition.section} we present the details of its detection edge.
Finally, our third aim is establishing asymptotic efficiency of the proposed estimator of the scale parameter and
intercept, when the intercept is not penalized.

Let, for $i=1, \ldots , n$, $x_i \in \R^p$ be a row vector of input variables and $y_i \in \R$ be a response variable.
The linear model is
$$ y_i = x_i \beta^* + \sigma^* \xi_i , \ i=1 , \ldots , n , $$
with $\beta^* \in \R^p$ an unknown (column-)vector of regression coefficients, $\sigma^* >0$ an unknown 
scale parameter - or noise level -, and $\{ \xi_i\}_{i=1}^n $ unobservable, i.i.d.\ noise variables with a given log-concave density $f$.
The number of variables  $p $ is allowed to be much larger than the number of
observations $n$. 
We will assume a sparsity condition on $\beta$, namely that it has not too many non-zero coefficients,
see Condition \ref{s.condition}.
We apply a transform of the minus-log-likelihood  and invoke an $\ell_1$-penalty.  The $\ell_1$-penalty on the regression 
coefficients $\beta\in \R^p$ is equal to
$\lambda \| \beta \|_1$ with $\lambda >0$ a tuning parameter and $\| \beta \|_1 := \sum_{j=1}^p | \beta_j |$ the
$\ell_1$-norm of the vector $\beta$.

Let $l(y) := - \log f(y)$, $y \in \R$. At $(\beta , \sigma) \in \R^p \times \R_+$, the  minus-log-likelihood, scaled by $1/n$,  is 
$${R}_n (\beta, \sigma) :={1 \over n}  \sum_{i=1 }^n \ell_{\beta , \sigma} (x_i,y_i) , $$
where 
$$ \ell_{\beta, \sigma} (x,y) = l\biggl ( { y - x \beta \over \sigma } \biggr ) + \log \sigma , \ (x,y) \in \R^p \times \R. $$
One calls $R_n (\beta , \sigma)$ the ``empirical risk" at $(\beta , \sigma)$.
For the case $\sigma^*$ known, the Lasso based on minus-log-likelihood loss is
\begin{equation} \label{Lasso.equation}
\min_{\beta \in \R^p } \biggl \{  {R}_n ( \beta , \sigma^*)  + \lambda \| \beta \|_1 / \sigma^* \biggr \} , 
\end{equation}
where $\lambda $ is a universal i.e., known, tuning parameter.
When one applies quadratic loss, i.e.\ $l(y)= y^2 $, $y \in \R$, this is known as the (classical) Lasso (\cite{tibs96}).
The theory for the  classical Lasso  is well-developed, see \cite{vandeG08},
\cite{brt09}, and the monographs \cite{koltch09a}, \cite{pbvdg11}, \cite{hastie2015statistical} and
\cite{giraud2021introduction}. The problem of the choice of the tuning parameter when $\sigma^*$ is unknown has
been also extensively studied. The paper \cite{belloni2011square}, 
\cite{sunzhang11}, introduced the square root Lasso for quadratic loss to
deal with unknown $\sigma^*$. Also theory for cross-validated Lasso is 
derived, see e.g.\ \cite{chetverikov2021cross}.

For the case $\sigma^*$ unknown, the idea is to transform the empirical risk $R_n$ using a given transformation
$\phi:\ \R \rightarrow \R$ such that the problem becomes ``pivotal", meaning that with the transformed
$R_n$ one can choose the tuning parameter independent of $\sigma^*$. We take $\phi$ as the exponential function $\phi (u):= \exp [u] $, $u \in \R$. This leads to what we call  the ``exp-Lasso"
\begin{equation}\label{exp-Lasso.equation}
 ( \hat \beta , \hat \sigma ) := \arg \min_{\beta \in \R^p , \sigma > 0 } 
\biggl \{ \exp \biggl [  R_n (\beta , \sigma) \biggr ] + \lambda \| \beta \|_1 \biggr \}. 
\end{equation} 
The choice $\phi = \exp [ \cdot ] $ allows to perform a disappearance act. Indeed, one may
write
$$ R_n ( \beta , \sigma) = \tilde R_n ( \tilde \beta , \tilde \sigma) + \log \sigma^* $$
where 
$ \tilde \beta = \beta / \sigma^* $, $\tilde \sigma = \sigma / \sigma^*$.
Thus
$$ \exp[R_n (\beta , \sigma) ] + \lambda \| \beta \|_1 = \sigma^* 
\biggl \{ \exp [\tilde R_n ( \tilde \beta , \tilde \sigma) ] + \lambda \| \tilde \beta \|_1 \biggr \} , $$
that is, in the theory, the minimization problem does not depend on $\sigma^*$. See Subsection
\ref{disappearance-act.section} for some more details.

A special case is Gaussian noise, where $f$ is the standard (say) normal density.
One easily verifies that in this special case, the exp-Lasso is the square root Lasso (see Subsection
\ref{Gaussian.example}). 
Our results are an extension to more general noise distributions. It is to be noted however that
we will require more sparsity than needed in the Gaussian case, see Condition \ref{s.condition}. 
This is due to our handling of the non-linearity of the problem for the non-Gaussian case. 
Examples include the Subbotin distribution (see \cite{olea2022out}), the logistic distribution, 
Huber's distribution, and the Gumbel distribution.
These examples will be treated in Section \ref{examples.section}, where
we also discuss the consequences when the noise distribution is misspecified. 

\subsection{Organization of the paper}
The main conditions and result can be found in  the next section (Section \ref{main.section}). In Section \ref{not-penalized.section}
we briefly discuss the adjustment when certain coefficients (e.g.\ the constant term) are not penalized.
In Section \ref{phase+.section} we examine variable selection and the detection edge.
We establish asymptotic efficiency of the estimator of the scale parameter and
the constant term in Section \ref{efficiency.section}.  Section \ref{examples.section} looks at examples, and in particular what
can be said in case of a misspecified noise distribution. 
Section \ref{proof-main.section} has the proof of the main result. Section \ref{further-proofs.section}
has the proofs of the results in Sections \ref{phase+.section} and
\ref{efficiency.section}. 

\subsection{Some notation}\label{notation.section}

The $\ell_1$-norm of a vector $b \in \R^p$ is $\| b \|_1 := \sum_{j=1}^p | b_j |$, its $\ell_2$-norm is
$\| b \|_2 := \sqrt {\sum_{j=1}^p b_j^2 }$ and its $\ell_{\infty}$-norm is $\max_{1 \le j \le p} | b_j| $.
For a matrix $A$, we write the maximum absolute value of is entries as $\| A\|_{\infty} $. 

We let $S^*:= \{ j \in \{ 1 , \ldots , p\} :\ \beta_j^* \not= 0 \} $ be the active set of $\beta^*$. For a vector
$ b \in \R^p$ we write $b_{S^*} := \{ b_j : \ j \in S^* \} $ and $b_{-S^*} := \{ b_j : \ j \notin S^* \} $.

In the proofs, we employ the following notation.
For a function $g: \R^{p+1} \rightarrow \R$, we write
$$ P_n g := {1 \over n} \sum_{i=1}^n g (x_i , \xi_i) , \ Pg := {1 \over n} \sum_{i=1}^n \EE g ( x_i, \xi_i) . $$
We apply the re-parametrization $b= (\beta- \beta^*)/ \sigma$ and $d= \sigma^* / \sigma$, $(\beta, \sigma) \in \R^p \times
\R_+ $. 
Thus, for $i=1 , \ldots , n$, 
$$ l\biggl ( { y_i - x_i \beta \over \sigma} \biggr ) = l ( d \xi_i - x_i b) =: g_{b,d} ( x_i , \xi_i ) , \ (b,d) \in \R^p \times \R_+ . $$
Write $\dot g_{b,d}^{\rm b}$ for the derivative of $g_{b,d}$ with respect to $b$ and
$\dot g_{b,d}^{\rm d} $ for its derivative with respect to $d$.

We let $M = M_n >0$ be a  sequence to be specified (see Theorem \ref{main2.theorem} for its full specification), tending to zero,  and define
$$ \Theta_M := \biggl \{ ( b , d ) \in \R^{p} \times \R_+: \ \| b \|_1 + | d-1 | \le M \biggr  \} , $$
where $\R_+ := (0, \infty)$.

\section{Main result}\label{main.section}

First, we state Conditions \ref{f.condition}, ..., \ref{s.condition}. In Theorem \ref{main.theorem} these
are used to derive an oracle inequality.

\begin{condition}\label{f.condition} The noise variables 
$\{ \xi_i \}_{i=1}^n $ are the first $n$ of an infinite sequence of i.i.d.\ copies of a
random variable $\xi$ with (known) density $f$. This density is strictly positive everywhere and
$l (\cdot) := - \log f (\cdot)$ is convex and differentiable with derivative $\dot l (\cdot)$. We furthermore assume that
$| \EE l (\xi)| < \infty $, $\EE ( \dot l (\xi) )^2 < \infty $ and $\EE ( \dot l (\xi) \xi )^2 < \infty $.
\end{condition}

Note that ${\rm var} ( \dot l (\xi) )= \EE ( \dot l (\xi) )^2$  is the Fisher information for location
and ${\rm var} ( \dot l (\xi) \xi ) = \EE ( \dot l (\xi) \xi )^2 -1$ is the Fisher information for scale.

\begin{condition} \label{Kx.condition}The co-variables  are fixed (i.e.\ non-random) and  bounded: for a constant $ K_{\rm x}\ge 1 $,
$$ \max_{1 \le i \le n} \max_{1 \le j \le p } | x_{i,j} | \le K_{\rm x} . $$
\end{condition}

One may argue that one can without loss of generality assume that the constant $K_{\rm x}$ is equal to one.
On the other hand, alternatively to fixed design, one may consider the situation with i.i.d.\  random design independent of $\{
\xi_i \}$. 
In the latter case, the co-variables are required to be bounded by some constant $K_{\rm x}$ with high probability.
We primarily have in mind here the case of Gaussian design. To avoid digressions we study this case only later, in Subsection \ref{random-design.section}. We remark
furthermore that Condition \ref{Sigma.condition} below is best understood in the context of random design.

The next condition is a local Lipschitz condition, locally near $\xi$,  on the derivative $\dot l$, with Lipschitz constant $G(\xi)$ depending
on the location $\xi$. For cases where it does not hold, we will discuss in Subsection \ref{notG.section} a similar result
as in Theorem \ref{main.theorem}, but with a (universal) tuning parameter that stays away from the detection edge. 

\begin{condition} \label{G.condition} For some function $G>0$ we have for $M$ small enough
$$|  \dot l (\xi + y) - \dot l (\xi + \tilde y ) | \le G(\xi) | y - \tilde y | , \forall  \ | y | \vee | \tilde y | \le(  K_{\rm x} + |\xi |)M, $$
where $\EE G^2 (\xi) < \infty$ and $\EE G^2 ( \xi) \xi^4 < \infty$.
\end{condition} 

\begin{condition} \label{H.condition} Let, for $(c,d) \in \R \times \R_+ $, the function $H( c,d)$ be defined as
$H(c,d)  := \EE l ( d \xi -c ) $. For $|c| + |d-1|$ small enough, its Hessian
$\ddot H (c,d)$ exists  and is continuous, and $\ddot H(0,1)$ is positive definite, with smallest eigenvalue $\kappa_0>0$.
\end{condition}

We will show in Subsection \ref{excess-risk.section} that Condition \ref{H.condition} holds when
the second derivative $\ddot l $ exists and is strictly positive. Condition \ref{H.condition} however only requires
the expected value to be twice differentiable. Since taking the expected value has a smoothing effect,
Condition \ref{H.condition} can also hold when $\ddot l$ does not exists, as 
in Example \ref{Huber.example}. 

We now list our conditions involving asymptotics. The high-dimensional model changes with
the number of observations $n$, but the density $f$ is kept fixed, not depending on $n$. Asymptotic
statements are for $n \rightarrow \infty$. 
With the notation $u \stackrel{\le}{\sim} v$, where $(u,v) = (u_n, v_n) $ is a sequence of strictly positive numbers, we mean that   $\limsup_{n \rightarrow \infty} u_n / v_n < \infty$ Similarly, $u \asymp v$ means
$u \stackrel{\le}{\sim} v$ and $v \stackrel{\le}{\sim} u $. With $u =u_n = o(1)$ we mean that
$\lim_{n \rightarrow \infty} u_n =0 $. Then $u \ll v $ or $v \gg u $ means $u/v = o(1)$.
For $u_n$ not necessarily positive, $u_n = {\mathcal O} (v_n)$ is another notation for $|u_n| \stackrel{\le}{\sim}  v_n $. 

\begin{condition}\label{p.condition}
The number of variables $p$ tends to infinity, and $\log p / n \rightarrow 0 $.
\end{condition}

The first part of Condition \ref{p.condition} is invoked because for $p$ remaining bounded the theory is of a different flavor.

Let
$$ \hat \Sigma := {1 \over n} \sum_{i=1}^n x_i^T x_i \in \R^{p \times p } $$
be the (normalized) Gram matrix. The first part of the next condition holds, if for all $n$, $\{ x_i \}_{i=1}^n$ are
$n$ realizations of a random variable ${\bf x} \in \R^p$ with sub-Gaussian entries with 
constant $K_{\rm x}$ and with $\EE {\bf x}^T {\bf x} = \Sigma$. The entries of $\Sigma$ should then
not grow with $n$.

\begin{condition} \label{Sigma.condition} For some matrix $\Sigma\in \R^{p \times p} $, it holds that

$$ \max_{(j,k ) \in \{ 1 , \ldots , p \}^2}  | \hat \Sigma_{j,k} - \Sigma_{j,k} | \stackrel{\le}{\sim} K_{\rm x}^2 \sqrt {\log p / n } . $$
Furthermore, $\Sigma$ has with smallest eigenvalue $\Lambda_{\rm x}^2 >0$.
\end{condition}

Recall the notation $\| b \|_{\infty} := \max_{1 \le j \le b } | b_j |$, $b \in \R^p$. Set
\begin{equation}\label{lambda*.equation}
 \lambda^*  :=  \biggl \| {1 \over n } \sum_{i=1}^n \dot l (\xi ) x_i \biggr \|_{\infty} \exp \biggl [ - {1 \over n} \sum_{i=1}^n \log f(\xi_i) \biggr ] , 
 \end{equation}
and let $F$ be the distribution of $\lambda^*$.

In the next condition, we either take take $0< \alpha <1/2$ fixed, not depending on $n$, or
(say) $\alpha = 1/p$. A fixed $\alpha$ is in line with our theory concerning the detection edge,
see Lemma \ref{phase-transition1.lemma}. The asymptotic confidence level in the result of Theorem \ref{main.theorem} will be $1-\alpha$. The choice $\alpha=1/p$ means that results hold with probability tending to
one. This is the right context for showing asymptotic efficiency of the estimator of $\sigma^*$. 

The next condition ensures that we can take the tuning parameter $\lambda 
\asymp  \sqrt {\log p / n}  $ but not of smaller order. This has to do with Condition \ref{Sigma.condition}:
the difference between the entries of $\hat \Sigma$ and $\Sigma$ is then not essentially larger than $\lambda$. 
 
\begin{condition} \label{alpha.condition} We have
$$ F^{-1} (1 - \alpha) \asymp  \sqrt {\log p / n} . $$
\end{condition}
See Lemma \ref{lambda*.lemma} for a justification of the upper bound in this condition.

Note that under Condition \ref{f.condition} $F^{-1} (1- \alpha)$ is a known constant.
It  will in general not be given in explicit form, but one can do
a Monte Carlo simulation to approximate it with any prescribed precision.
The tuning parameter $\lambda$ will be chosen larger than but asymptotically equal to $F^{-1} (1- \alpha)$.
We note that if the noise distribution is misspecified, and (partly) unknown, then $F$ is (partly) unknown
so that the problem of the choice of the tuning parameter is back again. Otherwise, our results do
not rely on a well specified noise distribution. See Subsection \ref{misspecified.section} and the examples
in Section \ref{examples.section} for some more
details.

\begin{condition} \label{s.condition} We assume that $ s^* \le s_{\rm max}$ where $s^* := \# \{ \beta_j^* \not= 0 \} $ and
where $1 \le s_{\rm max} \ll  \Lambda_{\rm x}^2  \sqrt {n / \log p}   /  K_{\rm x}^2$.  
\end{condition}

\begin{theorem} \label{main.theorem} Assume Conditions \ref{f.condition}, ..., \ref{s.condition}.
Let $0< \eta < 1$, $1- \eta \stackrel{\ge}{\sim} 1 $ and 
$$ \eta^2 \gg s_{\rm max} \sqrt {\log p / n}   K_{\rm x}^2  /\Lambda_{\rm x}^2    . $$ Take  
$\lambda \asymp \sqrt {\log p / n} $, $ \lambda   \ge F^{-1} (1- \alpha)/ (1- \eta) $.
Then with probability at least $1- \alpha + o(1)$,
\begin{eqnarray*}
{(\hat \beta - \beta^*)^T \hat \Sigma  (\hat \beta- \beta^*)\over \sigma^{*2}  }  &\stackrel{ \le} {\sim} & {      s^*  \lambda^2 \over  \Lambda_{\rm x}^2 }  + \lambda^2   ,
\end{eqnarray*}
and
\begin{eqnarray*}
{ \| \hat \beta - \beta^* \|_2^2 } &\stackrel{ \le } { \sim}   &  {\lambda^2 s^*  \over \Lambda_{\rm x}^2 }  + { \lambda^2 \over \Lambda_{\rm x} }     , 
\end{eqnarray*} 
as well as the (rough) bound
$$ { | \hat \sigma - \sigma^* |  \over \sigma^* } \stackrel{\le} {\sim} {\lambda \sqrt {s^*} \over \Lambda_{\rm x} } + \lambda ,$$
and finally also
$$ {\| \hat \beta - \beta^* \|_1 \over \sigma^*} \le M , $$
where $M = {\mathcal O}( \lambda s^* / ( \eta \Lambda_{\rm x}^2 ) + \lambda / \eta) $. 
\end{theorem}

If in Theorem \ref{main.theorem}, $\eta \rightarrow 0$, we say that we are (asymptotically)
at the detection edge, and if $\eta = 1- 1/C$ with $C>1$ a constant not depending on $n$,
we say we stay away from the detection edge. See Lemmas \ref{phase-transition1.lemma}
and \ref{phase-transition2.lemma} for the basis  of this way of saying.

\begin{remark} \label{s*=0.remark} A special case in Theorem \ref{main.theorem} is when there is no signal: $s^* =0$. 
The second term in the inequalities then starts playing its role.
\end{remark} 

%

\begin{remark} In Theorem \ref{main.theorem} we took either $0<\alpha<1/2$ not depending on $n$ or 
$\alpha = 1/p$. In the latter case the confidence level is at least $1- \alpha  +o(1) = 1- o(1)$ and we make no precise
statements how close the confidence level is to 100 \%.
A fixed value for $\alpha$ not depending on $n$ leads to a ``practical" choice for the tuning parameter.
\end{remark}

\section{Some coefficients not penalized}\label{not-penalized.section} 

There may be some input variables which are a priori considered as being necessarily included in the regression equation. For $x_j = ( x_{i,1} , \ldots , x_{n,j})^T$, $j=1 , \ldots , p$, let for $q < p$, $x_1 , \ldots , x_q$ be these 
necessary variables. There is no penalty on their coefficients. If we write
\begin{equation}\label{x.equation}
x_i =\underbrace{(x_{i,1} , \ldots, x_{i,q})}_{x_{i,0}} , \underbrace{( x_{i,q+1} , \ldots , x_{i,p} ) }_{x_{i, -0}} =: (x_{i,0} , x_{i, -0}) ,
\end{equation}
and 
\begin{equation}\label{beta.equation}
\beta^T = \underbrace{( \beta_1 , \ldots , \beta_q )}_{\beta_0^T}, 
\underbrace{( \beta_{q+1} , \ldots , \beta_p )}_{\beta_{-0}^T } =: ( \beta_0^T , \beta_{-0}^T) , 
\end{equation}
the new exp-Lasso is
\begin{equation}\label{new-exp-Lasso.equation}
 (\hat \beta  , \hat \sigma) := 
\arg \min_{\beta \in \R^{p} , \ \sigma >0 } \biggl \{ 
\exp[ R_n ( \beta, \sigma) ] + \lambda \| \beta_{-0} \|_1 \biggr \} .  
\end{equation}

A lazy way to deal with this new exp-Lasso is
by viewing $(x_1 , \ldots , x_q)$ as active variables,  thus obtaining a newly defined active set 
\begin{equation}\label{new-active-set.equation}
S_+^* := \{ 1, \ldots , q\} \cup \{ \beta_j^* \not= 0 , \ q+1 \le j \le p \}.
\end{equation}
One obtains the following corollary of Theorem \ref{main.theorem}.

\begin{corollary}  Suppose the conditions of Theorem \ref{main.theorem} are met with the newly
defined active set $S_+^*$ given in (\ref{new-active-set.equation}), and with $s^*$ 
replaced by $s_+^* := | S_+^* | $, the cardinality of this new 
active set. Then the conclusion of Theorem \ref{main.theorem} is valid for the exp-Lasso given
in (\ref{new-exp-Lasso.equation}). 
\end{corollary}

\section{Variable selection and the  detection edge}\label{phase+.section}

From now on, we assume for (simplicity) that $\Lambda_{\rm x}$ in Condition \ref{Sigma.condition} does not depend on $n$. We also assume that $\sigma^*=1$, which can be
done without loss of generality by the disappearance act.

We study in this section the new exp-Lasso, where $x_{i,1} =1$ for all $1 \le i \le n $, i.e., we include an
intercept which is not penalized in the regression equation. We can then take 
each $x_j$ with $j \ge 2$ in deviation from its mean $\bar x_j = \sum_{i=1}^n x_{i,j} /n $, $j=2 , \ldots , p $.
Instead of changing the notation, we assume that $\bar x_j = 0$ for
$j=2 , \ldots , p$. We let
$$x_i =\underbrace{(1) }_{x_{i,0}} , \underbrace{( x_{i,2} , \ldots , x_{i,p} ) }_{x_{i, -0}} =:(x_{i,0} , x_{i, -0}) $$
and 
$$\beta^T = \underbrace{ (\beta_1 ) }_{\beta_0}, \underbrace{( \beta_{2} , \ldots , \beta_p )}_{\beta_{-0}^T} =: ( \beta_0 , \beta_{-0}^T) . $$


We write ${\cal H} (c,d) := H(c,d) - \log d$, $(c,d)\in \R\times \R_+ $, and, whenever the second
derivative exist, 
$$ \ddot {\cal H}(c,d):= \begin{pmatrix} \ddot{\cal H}^{{\rm d} , {\rm d} } (c,d) &  \ddot{\cal H}^{{\rm c}, {\rm d}}  (c,d) \cr 
\ddot{\cal H}^{{\rm c} , {\rm d}} (c,d) & \ddot {\cal H}^{{\rm c}, {\rm c}} (c, d) \cr \end{pmatrix} , $$
where 
$\ddot{\cal H}^{{\rm d} , {\rm d} } $ is the second derivative with respect to $d$, $\ddot {\cal H}^{{\rm c}, {\rm c}} $
the second derivative with respect to $c$ and $\ddot{\cal H}^{{\rm c} , {\rm d}} $ the mixed derivative.

\subsection{Variable selection under the irrepresentable condition}

\begin{condition}\label{H2.condition}
For some constant $L_{\rm H}$ and for all $|c| + | d-1|$ small enough,
$$ \| \ddot {\cal H} (c,d) - \ddot {\cal H} (0,1)\|_{\infty} \le L_{\rm H} (| c| + |d-1| ) . $$
\end{condition}

  Let $S^* := \{ j \ge 2: \ \beta_j^* \not= 0 \} $ be the active set of the penalized coefficients and
  $s^* := | S^* | $.
 The matrix $X_{S^*} $ is defined as selecting only the
 columns in $S^*$ and $X_{-S}^* $ selects only the columns in $\{ 2 , \ldots , p \} \backslash S^* $.
 
 The irrepresentable condition was introduced in \cite{ZY07}  (see also \cite{mebu06}) and used
 for variable selection with the classical Lasso which is based on quadratic loss. 
 
 \begin{condition}\label{irrepresentable.condition} The matrix $X_{S^*}^T X_{S^* } $ is invertible.
 For some $0 \le \eta_0 \le1$ and for all vectors $\tau_{S^*} \in
 \R^{s^*}  $ with $\| \tau_{S^*} \|_{\infty} \le 1$,  it holds that
 $$ \| X_{-S^*}^T X_{S^*} ( X_{S^*}^T X_{S^* } )^{-1}  X_{S^*} \tau_{S^*} \|_{\infty} < \eta_0 . $$
 \end{condition}
 
   \begin{theorem} \label{variable-selection.theorem} 
   Suppose the conditions of Theorem \ref{main.theorem} are met
   with $\Lambda_{\rm x}$ not depending in $n$. Take  $\sum_{i=1}^n x_{i,-0} =0$. Assume
   Conditions \ref{H2.condition} and \ref{irrepresentable.condition} as well, with 
   $$\eta_0= { \eta(1- r_n) \over 2 - \eta (1+r_n)  }, $$
   where $0< r_n = {\mathcal O} ( K_{\rm x}^2 \lambda s^* / \eta^2) = o(1) $.
   Then $\hat \beta_{-S^* } = 0 $  with probability at least $1- \alpha + o(1)$.
   \end{theorem}
 
    Note that the above result favors a value of $\eta$ close to 1. On the other hand, a value of $\eta$ close to
   zero allows the tuning parameter $\lambda$ to be close to the detection edge.  
   
   \subsection{The detection edge}{\label{phase-transition.section}

\begin{lemma}\label{phase-transition1.lemma} Assume Conditions \ref{f.condition}, ..., \ref{s.condition},
with
$\Lambda_{\rm x}^2 $ fixed and $\sum_{i=1}^n x_{i,-0} =0 $. 
Under $H_0 : \ \beta_{-0}^* =0 $, it holds that $\hat \beta_{-0} =0$
with probability at least $1-\alpha + o(1)$.
 \end{lemma}

\begin{lemma}\label{phase-transition2.lemma}
Assume Conditions \ref{f.condition}, ..., \ref{s.condition}, with $\Lambda_{\rm x}$ fixed. We also
require  Condition \ref{H2.condition} and  that $\sum_{i=1}^n x_{i,-0} = 0$ and 
$\min_{2 \le j \le p } \| x_j \|_2^2 /n  \stackrel{\ge}{\sim} 1  $. Let $H_0: \beta_{-0}^*=0$ be true.
Then for $1 > \eta \gg K_{\rm x}^2 \sqrt {\log p / n} $,
$$ \PP ( \hat \beta_{-0} \not= 0 ) \ge \PP ( \lambda^* (1- \eta) > \lambda ) + o(1) . $$
\end{lemma}

\section{ Asymptotic efficiency}
\label{efficiency.section}

In this section we assume throughout, and without loss of generality, that $\sigma^*=1$. 
When applying Lemma \ref{efficiency.lemma} e.g.\ for building asymptotic confidence intervals for
the scale parameter $\sigma^*$ and the intercept $\beta_0^*$ one then should use the proper rescaling. 

\subsection{Scale parameter and intercept not penalized}\label{scale+intercept.section}

We apply the new exp-Lasso, where $x_{i,1} =1$ for all $1 \le i \le n $, i.e., we include an
intercept. We assume the intercept is not penalized in the regression equation. 

We let, for $(m, d) \in \R\times \R_+$,   $ {\cal K} (m,d) :=
\EE \log f ( \xi ) - (\EE \log f (d (\xi-m)) - \log d )$ be the Kullback-Leibler information.

Suppose Condition \ref{H2.condition} holds. Let 
$$
 \ddot {\cal K} (m,d) = \begin{pmatrix} \ddot {\cal K}^{{\rm d}, {\rm d} } (m,d) & \ddot {\cal K}^{{\rm d}, {\rm m} } (m,d) \cr
\ddot {\cal K}^{{\rm d}, {\rm m} } (m,d) & \ddot {\cal K}^{{\rm m}, {\rm m} } (m,d)  \cr \end{pmatrix} ,$$
be the Hessian of ${\cal K} $ at $(m,d)$ (whenever it exists).
Then $\ddot {\cal K} (0,1)= \ddot {\cal H} (0,1)$. 
Under Condition \ref{H2.condition}, it is true that 
for some constant $L_{\rm K}$  and for all
sufficiently small $| m| + |d-1| $ that
$$ \| \ddot {\cal K} ( m,d)  - \ddot {\cal K} (0,1)  \|_{\infty} \le L_{\rm K} (| m| + |d-1|) . $$

%
%
%
%

The next lemma shows that under a slightly stronger condition on $s_{\rm max}$, the exp-Lasso estimator of 
the scale parameter $\sigma^*$ and
the intercept $\beta_0^*$ are asymptotically equal to the MLE of these parameters 
when $\beta_{-0}$ were known, so that
$$\sqrt n \begin{pmatrix}  \hat d -1 \cr \hat \beta_0 - \beta_0^* \cr \end{pmatrix} 
\stackrel{{\cal D}}{\rightarrow} {\cal N} (0, \ddot {\cal K}^{-1} (0,1)) , $$
where $\stackrel{\cal D}{ \rightarrow}$ means convergence in distribution and
where ${\cal N} (0, \ddot {\cal K}^{-1} (0,1)) $ is the 2-dimensional normal distribution
with mean zero an co-variance matrix $\ddot {\cal K}^{-1} (0,1)$.

\begin{lemma}\label{efficiency.lemma} Assume
Condition \ref{H2.condition}.  Suppose moreover the conditions of Theorem \ref{main.theorem}  hold
with $\alpha = 1/p$, with $\Lambda_{\rm x}^2$ not depending on $n$ and 
where
Condition \ref{s.condition} is strengthened to 
$$s_{\rm max}   \ll \min\{ \sqrt n / \log p , \sqrt { n/ \log p} / K_{\rm x}^2 \} .  $$ 
Let $\sum_{i=1}^n x_{i, -0} =0$.
Then we have with probability tending to one, 
$$ \begin{pmatrix}  \hat d -1 \cr \hat \beta_0 - \beta_0^* \cr \end{pmatrix}=- 
\ddot {\cal K}^{-1}(0,1)    
P_n \begin{pmatrix}  { \partial  \ell_{\beta, \sigma}  \over \partial d} \biggl \vert_{\beta= \beta^*,  d= 1}\cr
 { \partial  \ell_{\beta, \sigma}  \over \partial \beta_0 } \biggl \vert_{\beta= \beta^* , d= 1}\cr
\end{pmatrix}  
+ o(n^{-1/2} )  .$$
\end{lemma}

  \subsection{Further not-penalized parameters}\label{efficiency-general.section}
  More general than in the previous subsection, take
  the exp-Lasso as
\begin{equation}\label{new-exp-Lasso2.equation}
 (\hat \beta  , \hat \sigma) := 
\arg \min_{\beta \in \R^{p} , \ \sigma >0 } \biggl \{ 
\exp[ R_n ( \beta, \sigma) ] + \lambda \| \beta_{-0} \|_1\biggr \}  .  
\end{equation}
with for $i=1 , \ldots , n $,
$x_i =(x_{i,0} , x_{i, -0}) $ defined as in equation (\ref{x.equation}), and with 
$\beta = ( \beta_0^T , \beta_{-0}^T) $ defined as in equation (\ref{beta.equation}).
We assume $q$ is fixed, not depending on $n$.
Write
$$ X = \begin{pmatrix}  x_{1} \cr \vdots \cr  x_{n} \cr \end{pmatrix} =: ( X_0 , X_{-0}) \in \R^{n \times p }, $$
Let $\Pi$ be the projection operator on the space spanned by the columns of
$$X_{0}:= \begin{pmatrix}  x_{1,0} \cr \vdots \cr  x_{n,0} \cr \end{pmatrix} \in \R^{n \times q } . $$
Then
$$ X\beta^* = X_0 \beta_0^* + X_{-0} \beta_{-0}^* = X_0 \beta_0^* + \Pi X_{-0} \beta_{-0}^* +
(I- \Pi ) \beta_{-0}^* . $$
$$ = X_0 \gamma_0^* + (I- \Pi ) \beta_{-0}^* , $$
where
$$ \Pi X_{-0} = X_0 \Gamma, \ \Gamma \in \R^{q \times (p-q)  } , \ \gamma_0^* = \beta_0^* + \Gamma \beta_{-0}^* . $$
In other words, under the conditions of Lemma \ref{efficiency.lemma}, the estimator of $\gamma_0^*$ can be shown to be asymptotically efficient using the
same arguments as in the previous subsection. Nonetheless, unless $\Pi X_{-0} =0 $, the new parameter
$\gamma_0^*$ may in practice not be the parameter of interest.

 \section{Examples}\label{examples.section}
 
 Before looking at examples, it may be relevant to discuss what can be said
 when the conditions of Theorem \ref{main.theorem} are not satisfied.
 First we briefly look at random design.
 
 \subsection{Random design} \label{random-design.section} Suppose that $\{ x_i \}_{i=1}^n $ are $n$ realizations of a random row-vector ${\bf x} \in \R^p$. Then Condition \ref{Kx.condition} holds with $K_{\rm x} =1$ if $\| {\bf x } \|_{\infty} \le 1$ and
Condition \ref{Sigma.condition} holds with probability tending to one if $\Sigma := \EE {\bf x}^T {\bf x} $ has smallest eigenvalue $\Lambda_{\rm x}^2$. 
Alternatively, when ${\bf x}$ is a standard (say) Gaussian random vector, then Condition
\ref{Kx.condition}, with $K_{\rm x} = {\mathcal O} ( \sqrt {\log (np) })$, as well as Condition
\ref{Sigma.condition}, with $\Lambda_{\rm x} =1$, are met with probability tending to one.
We then require in Condition \ref{s.condition} that
$$ s^* \sqrt {\log p / n } \log (np)  = o(1) .$$ 
Moreover, in this case  Condition \ref{irrepresentable.condition} holds with probability tending
to one for $\eta_0= {\mathcal O} ( s^*\sqrt {\log p / n} ) $. Thus, not surprisingly, with i.i.d.\
standard Gaussian design, one can get near the detection edge and (yet), by
Theorem \ref{variable-selection.theorem},  do variable selection.

In our examples, we mainly look at Condition \ref{f.condition} which assumes $f$ is the true density of the noise,
and Condition \ref{G.condition} which assumes  $\dot l $ is Lipschitz with appropriate Lipschitz constant depending
on location. 
 
  \subsection{Misspecified noise distribution} \label{misspecified.section}
  Suppose that $f$ is possibly not the density of the noise.
  Let $ f^*$ be the density of $\xi$ and suppose $f^*$ is in part unknown. 
  Write $\EE_{f^*}$ for expectation under the distribution with density $f^*$.
The theory goes through  if $f\not= f^*$ under moment conditions on $f^*$. The
problem is however that the choice of the tuning parameter as $(1- \alpha)$-quantile of the
distribution of $\lambda^*$ is no longer possible, as it depends on the unknown
distribution of $\xi$. 
We however have the normalization
$$ \EE_{f^*} \dot l (\xi) = \int \dot l (y ) f^*(y) =0 , \ \EE_{f^*} \dot l (\xi) \xi = \int \dot l (y) y f^* (y)=1 ,$$
i.e., we do know something about $f^*$. 
For the choice of the tuning parameter $\lambda$ we need upper bounds for
$$ \EE_{f^*} ( - \log f (\xi ) )  \  {\rm and} \ \EE_{f^*} (\dot l (\xi))^2 . $$
An upper bound for the first expection is
$$  \EE_{f^*} ( - \log f (\xi ) ) \le \min_{y} l(y) + \underbrace{\int \dot l (y) y f^* (y) }_{=1} . $$
If also an upper bound for the second expectation is available, we call the model ``robust".

\subsection{Condition \ref{G.condition} violated} \label{notG.section}
Condition \ref{G.condition} is a Lipschitz condition on the derivative of $l= - \log f $ with Lipschitz constant
depending on location. If it is not true one probably has to let go the ambition to have a choice
of the tuning parameter $\lambda$ close to the detection edge do its job. However, one may still
have good results when one takes $\lambda$ larger.

The following condition ensures that the result of Theorem \ref{main.theorem} remains true,
say for $\alpha=1/p$ and the condition on $\lambda$ given by
$\sqrt {\log p / n} \asymp \lambda \ge C \sqrt {\log p / n}  $ (i.e.\ $\eta \ge 1- 1/ C$), where $C $ is a known fixed
(not depending on $n$)
constant.

\begin{condition} \label{G0.condition} For some function $G_0>0$ we have for $M$ small enough
$$|  l (\xi + y) - l (\xi + \tilde y ) | \le G_0(\xi) | y - \tilde y | , \forall  \ | y | \vee | \tilde y | \le(  K_{\rm x} + |\xi |)M  $$
where $\EE G_0^2 (\xi) < \infty$, and $\EE G_0^2 (\xi) \xi^2 < \infty$.
\end{condition}

 \subsubsection{Gaussian noise distribution} \label{Gaussian.example}
 In this case, by straightforward manipulation,
 $$ \min_{\beta\in \R } \biggl \{ \min_{\sigma >0 } \exp[ R_n ( \beta , \sigma) ] \biggr \} + \lambda \| \beta \|_1 $$
 $$ = {\rm e}^{(1+ \log (2
 \pi) )/2 } \min_{\beta \in \R}  \left ( \sum_{i=1}^n ( y_i - x_i \beta )^2 / n \right )^{1/2}+ \lambda \| \beta\|_1 . $$
 In other words, the exp-Lasso is the square root Lasso. Of course, we can add any constant
 to the log-likelihood, i.e., the term
 ${\rm e}^{(1+ \log (2
 \pi) )/2 }$ can be neglected. 
 
 We argue that if the noise distribution is misspecified, but with
 finite first and second moment, the misspecification has
 almost no impact especially when the noise is also symmetric. 
 Note that by the normalization of Subsection \ref{misspecified.section}, the distribution with density $f^*$ has
 $$ \EE_{f^*} \dot l ( \xi) = \EE_{f^*} \xi \stackrel{\triangle}{=} 0 , \ \EE_{f^*} \dot l(\xi ) \xi = \EE_{f^*} \xi^2 
 \stackrel{\triangle}{=}1 . $$
 So also
 $$ \EE_{f^*} (- \log f(\xi) ) = (1+ \log (2 \pi))/2 ,  \ \EE_{f^*} ( \dot l (\xi))^2 = \EE \xi^2 =1. $$
 Thus the model is ``robust" in the sense of Subsection \ref{misspecified.section}. 
 We only  have to avoid a slightly too optimistic choice for
 the tuning parameter. There is however a good universal bound (see Lemma \ref{lambda*.lemma}), at least for
 bounded fixed design or the random design as described in Subsection \ref{random-design.section}.
 Possibly one stays a bit away from the detection edge.  
  
  \subsubsection{Subbotin noise distribution} \label{Subbotin.example}
  The density of the standard Subbotin distribution is
  $$ f(y) ={ r \over 2 r^{-1/r}  \Gamma (1/r) } \exp[- | y |^r/ r ] , \ y >0 , $$
  where $r >0$ is the shape parameter. We assume $r$ fixed, and $r\ge 1$ so that
  $y \mapsto - \log f(y) $ is convex. Now, ignoring the normalizing constant
  ${ r \over 2 r^{-1/r}\Gamma (1/r) }$, we get 
  $$ \EE_{f^*} \dot l (\xi) \xi =
  \EE_{f^*}  | \xi |^r \stackrel{\triangle}{=} 1 \ \Rightarrow  \  \EE_{f^*} ( - \log f(\xi) )= {1 / r} . $$
  Moreover
  $$ \EE_{f^*} (\dot l (\xi))^2= \EE_{f^*} \xi^{2(r-1)} . $$
  For a well-specified model $\EE_{f} \xi^{2(r-1)}$ is known, and if the model is not well-specified we
  have the bound $\EE_{f^*} \xi^{2(r-1)} \le 1$ for $1 \le r \le 2 $. For $1 \le r < 2$ however, Condition
  \ref{G.condition} does not hold. We replace it by Condition \ref{G0.condition}, where for $r>1$ we can take
  $$G_0 (\xi) \asymp  | \xi |^{r-1}\{ |\xi | > M^*\}  + M^{*(r-1) } \{ | \xi | \le M^* \} $$
  with $M^* = {\mathcal O} (  \sqrt { \log p / n }  s^* / ( \eta \Lambda_{\rm x}^2 )) $.
  Thus, for $1< r \le 2$ we have ``robustness" as for the case $r=2$, but we do not get
  near the detection edge.
  For $r=1$ Condition \ref{H.condition} is not satisfied.
 However, the estimator for general $r\ge 1$ is (avoiding constants in the likelihood)
 $$ \hat \beta = \arg \min_{\beta \in \R }
 \left ( \sum_{i=1}^n ( y_i - x_i \beta )^r / n\right )^{1/r} + \lambda \| \beta \|_1 , $$
 which is convex problem. This means Condition \ref{H.condition} is not necessary for the Subbotin case:
 one can replace $H(c, d) = \EE_{f^*} l( d \xi - c)$ by ${\cal H}(c,d) = \EE_{f^*} l( d \xi - c) - \log d $,
 $(c,d) \in \R \times R_+$. \\
 The reason for taking a Subbotin error distribution with $1 \le r <2$ may indeed be that one aims at robustness, without actually believing that this distribution is the true error distribution.

 \subsubsection{Logistic noise distribution}\label{logistic.example}
 
 The logistic distribution has density
$$ f(y) = { {\rm e}^{-y} \over (1+ {\rm e}^{-y} )^2 } , \ y \in \R . $$. 
Thus, $f$ is symmetric, 
$$ l(y) := - \log f (y) = y +2 \log (1+ {\rm e}^{-y} )  $$
and 
\begin{eqnarray*}
 \dot l(y) = 1- { 2 {\rm e}^{-y} \over 1+ {\rm e}^{-y} } =  {2 \over 1+ {\rm e}^{-y} } -1,\ 
 \ddot l (y) = { 2 {\rm e}^{-y} \over (1 + {\rm e}^{-y} )^2 } \ge 0 , \ y \in \R . 
 \end{eqnarray*}
 We see that $l$ is convex, that
 $\| \dot l \|_{\infty} \le 1 $ and  that $\| \ddot l \|_{\infty} \le 1/2$. 
 It follows that Condition \ref{G.condition} is satisfied with $G(\cdot )\equiv 1/2 $, provided
 $\EE_{f*} \xi^4 < \infty$.
 As we have $\EE_{f^*} (\dot l ( \xi))^2 \le 1$, we conclude that the model is ``robust".
 
  \subsubsection{Huber noise distribution} \label{Huber.example} 
  Here we take
 $$ l(y) = \begin{cases}  y^2/2 ,  & | y | \le 1 \cr
  (|y|  - 1/2) ,  & |y| \ge 1 \cr 
 \end{cases}, \ \ \ y \in \R.
 $$
 Then
 $$ \dot l (y) =  \begin{cases}  y , & | y | \le 1 \cr + 1 , & y \ge 1 \cr -1 , & y \le 1 \cr \end{cases}, \ \ \ y \in \R . $$
 Since $\dot l (\cdot)$ is Lipschitz with Lipschitz constant 1, one  sees that in Condition \ref{G.condition},
 one can take $G(\cdot) \equiv 1$. 
 We see that
 $$ \EE_{f^*} \dot l (\xi) \xi = \EE_{f^*} y^2 \{ | y | \le 1\}  + \EE_{f^*}  |y| \{ | y | >1 \} \stackrel{\triangle}{=} 1. $$
 Thus
 $$ \EE_{f^*}  (-\log f (\xi) ) \le 1/2 . $$
 Also $\EE_{f^*} (\dot l(\xi))^2 \le 1$. So the model is ``robust". 
 Note that 
 $$ \lim_{\sigma \downarrow 0} \sigma l (y/ \sigma) = | y | , \ y \in \R . $$
 Our context is different as for location parameter
 $m\in \R$ and scale parameter $\sigma>0$, we are looking at 
 $$ l(y-m/ \sigma) + \log \sigma , \ y \in \R .$$
 In our context, the point where the loss function goes from quadratic to linear
 is estimated.

 \subsubsection{Gumbel noise distribution}\label{Gumbel.example}
 The Gumbel distribution is used to model the distribution of extreme values.
 In this case
 $$ l(y) = y + {\rm e}^{-y} , \ 
  \dot l (y) = 1- {\rm e}^{-y} , \ 
  \ddot l (y) = {\rm e}^{-y} > 0 , \ y \in \R . $$
 So $l$ is convex.
  Condition \ref{G.condition} holds with 
  $$ G(\xi) = {\rm e}^{ - \xi (1+M) } (1+ {\rm e}^{K_{\rm x} M} ), $$
  provided $\EE_{f^*} G^(\xi) < \infty$ and $\EE_{f^*}  G^2 (\xi) \xi^4 < \infty $. The latter two
  moment conditions are met if the model is well-specified.

Note that for the well-specified case
$$ \EE_{f} (-\log f (\xi ))  = \gamma +1 $$
where $\gamma \approx 0.5772 $ is the Euler-Mascheroni constant.

Next, we discuss what happens when the noise distribution
is misspecified. The normalization of Section \ref{misspecified.section} says
$$ \EE_{f^*}{\rm e}^{\xi } \stackrel{\triangle}{=} 1 , \ \EE_{f^*} (1-{\rm e}^{-\xi} ) \xi \stackrel{\triangle}{=} 1 . $$
This implies
$$ \EE_{f^*} (- \log f(\xi) ) \le 2. $$
Moreover
$$ \EE_{f^*}( \dot l (\xi))^2 = \EE_{f^*} {\rm e}^{-2\xi} - 1 . $$
For the well-specified case $\EE_f {\rm e}^{-2\xi} =2 $.
We conclude that the model is not ``robust" in the sense of Subsection
\ref{misspecified.section}. Nevertheless, if one chooses the Gumbel
distribution as noise distribution, there generally is a reason for that.
If we know that we are not too far away from the Gumbel, say that for a
given $\epsilon>0$
$$ \EE_{f^*} {\rm e}^{-2\xi} \le 2+ \epsilon $$
we can use this in the choice of the tuning parameter to make the exp-Lasso robust
in this ``$\epsilon$-environment". This $\epsilon$-environment can be seen as quantifying that we know
that extreme negative values of the noise are rare.

\section{Proof of Theorem \ref{main.theorem}.} \label{proof-main.section} 

In this section we assume throughout Conditions \ref{f.condition}, ..., \ref{s.condition}. The only 
exception is Lemma \ref{H2.lemma} where we prove Conditions \ref{G.condition} and \ref{H.condition},
instead of assuming these. 

\subsection{The disappearance act}\label{disappearance-act.section}
The disappearance act is closely related to the concept of equivariance.
With the new notation we have for $(b,d) \in \R^p \times \R_+$, 
$$ R_n ( \beta , \sigma) = P_n g_{b,d} - \log d + \log \sigma^*, $$
with $\beta = \beta^*+\sigma^* b/d $ and $d= \sigma^* /\sigma $. Thus
$$ \exp[R_n ( \beta , \sigma]  + \lambda \| \beta \|_1 =
{\sigma^*} \biggl \{ \exp[P_n g_{b,d} - \log d ] + \lambda  \| \beta^* / \sigma^*  + b/d   \|_1 \biggr \} .$$
In other words
$$ (\hat b , \hat d) = \arg \min_{b\in \R, \ d>0 } \biggl \{ \exp[P_n g_{b,d} - \log d ] + \lambda  \| \beta^* / \sigma^*  + b/d   \|_1 \biggr \} . $$

\subsection{A basic inequality} \label{basic-inequality.section} 
Recall the empirical risk
$${R}_n (\beta, \sigma) :={1 \over n}  \sum_{i=1 }^n \ell_{\beta , \sigma} (x_i,y_i) , \  (\beta  , \sigma) \in \R^p \times \R_+ .$$
To establish convergence of a penalized empirical risk minimizer to the true value, one typically uses that
the penalized empirical risk at the estimator is smaller than or equal to the penalized empirical risk at the true value. This we call the basic
inequality. If the penalized empirical risk is convex in its parameters, this can be exploited to localize the problem,
that is, to get into an
appropriate neighborhood of the true value.
However, in our case $R_n ( \beta , \sigma)$
is not convex in $(\beta, \sigma)$. We can make it convex using another parametrization. We choose here
$b:= ( \beta- \beta^*) / \sigma$ and $d:= \sigma^* / \sigma $ and let $b^*=0$ and $d^*=1$. With this new parametrization nonetheless, 
the penalty $\lambda \| \beta \|_1 = \lambda \| \beta^* + \sigma^* b / d \|_1 $ becomes non-convex.
It turns out that, apart from being pivotal, the  $\phi$-transform $\phi [ \cdot ] = \exp [ \cdot ] $ of the  
empirical risk deals with the non-convexity.

Of course the parametrization $(\beta, \sigma) \mapsto (b,d)$ cannot be used for computing the exp-Lasso
$(\hat \beta , \hat \sigma)$. The computational problem generally remains non-convex.
 There are exceptions, the non-convexity problem disappears  for example when for a fixed $\beta$ the solution for the estimator $\hat \sigma_{\beta}$ of the exp-Lasso
is a convex function of $\beta$, as is the case for the square-root Lasso (and more generally
for Subbotin distributions). 

To turn the basic inequality into one convex in $(b,d)$ we use the following lemma.

\begin{lemma} \label{surrogate.lemma} We have for all scalars $u$ and $v$
$$ {\rm e}^{ v + \log \sigma} - {\rm e}^{u + \log \sigma^* }  \ge
\sigma^* {\rm e}^u  \biggl \{ {v-u +1 -d \over d}\biggr \}.$$
\end{lemma}

{\bf Proof.}  We first note that ${\rm e}^v- {\rm e}^u \ge {\rm e}^u (v-u) $. Thus
\begin{eqnarray*}
 {\rm e}^{v + \log \sigma} - {\rm e}^{u + \log \sigma^* }&=&
\sigma {\rm e}^v -\sigma^* {\rm e}^u  
 \stackrel{d= \sigma^* / \sigma}{=} { \sigma^*\over d}  {\rm e}^v - \sigma^* {\rm e}^u \\
& =& \sigma^* \biggl \{ {   {\rm e}^v  \over d} - {\rm e}^u \biggr \} 
 =  \sigma^* \biggl \{  { {\rm e}^v - {\rm e}^u   \over d}   + \biggl ({1 \over d} -1 \biggr ) {\rm e}^u \biggr \}  \\
&\stackrel{{\rm e}^v - {\rm e}^u \ge {\rm e}^u (v-u) }{ \ge } &\sigma^* \biggl \{ { {\rm e}^u (v-u) \over d}  +  \biggl ({1 \over d} -1 \biggr ) {\rm e}^u  \biggr \} \\
& = & \sigma^* {\rm e}^u  \biggl \{ {v-u +1 -d \over d}   \biggr \}. 
\end{eqnarray*}
\hfill $\sqcup \mkern -12mu \sqcap$

Recall for $b \in \R^p$ the notation $b_{S^*} := \{ b_j : \ j \in S^* \} $ and $b_{-S^*} = \{ b_j: \ j \notin S^* \} $.

\begin{theorem} \label{twisted.theorem} We have for all $0 \le t \le 1 $, and for  $\hat b_t ;= t \hat b + (1- t) b^* (= t \hat b ) $
and $\hat d_t :=  t \hat d + (1-t) d^* (= t \hat d + 1-t )$, the twisted basic inequality
\begin{equation}\label{twisted.equation}
 {\rm e}^{P_n g_{0,1} } \biggl \{P_n (g_{\hat b_t , \hat d_t} - g_{0,1} ) +1 - \hat d_t \biggr \}  \le 
\lambda ( \| \hat b_{t, S^*} \|_1   -\| \hat b_{t, -S^*} \|_1 )  .  
\end{equation} 
\end{theorem}

{\bf Proof of Theorem \ref{twisted.theorem}.}
First note that $R_n (\beta , \sigma) = P_n g_{b,d} + \log \sigma $. Therefore the basic 
basic inequality inequality for the exp-Lasso is
$$ {\rm e}^{ P_n g_{\hat b, \hat d} + \log \hat \sigma } - {\rm e}^{ P_n g_{0, 1} + \log \sigma^* }
\le   \lambda \| \beta^* \|_1 - \lambda \| \hat \beta \|_1. $$
On the other hand, by Lemma \ref{surrogate.lemma},
$${\rm e}^{P_n g_{\hat b, \hat d} + \log \hat \sigma } - {\rm e}^{P_n g_{0, 1} + \log  \sigma^* } 
 \ge \sigma^* {\rm e}^{P_n g_{0,1}} \biggl \{ { P_n (g_{\hat b , \hat \sigma} -  g_{0,1} )+
1 - \hat d  \over \hat d}   \biggr \}  ,$$
so that
$$ \sigma^* {\rm e}^{ P_n g_{0,1}} \biggl \{ { P_n (g_{\hat b , \hat d } -  g_{0,1} )
+1 - \hat d  \over \hat d}  \biggr \} \le \lambda ( \| \beta^* \|_1 - \| \hat \beta \|_1)  $$
or
$${\rm e}^{ P_n g_{0,1}} 
\biggl \{ P_n (g_{\hat b , \hat d } -  g_{0,1}) 
 + 1 - \hat d  \biggr \}   \le  { \lambda \hat d \over  \sigma^*} ( \| \beta^* \|_1 - \| \hat \beta \|_1) . $$
 Now recall that for $(\beta , \sigma) \in \R^p \times \R_+$ the re-parametrization $b = d (\beta - \beta^* ) / \sigma^*$ 
 (and 
 $d= \sigma^* / \sigma$) so that $\beta = \beta^* + b \sigma^* / d $. 
It follows that $\beta- \beta^*= b  \sigma^* / d $.
Thus
\begin{eqnarray*}
 \| \beta^* \|_1 -\| \beta \|_1 &= &\| \beta^*\|_1 - \| \beta_{S^* }\|_1 - \| \beta_{-S^* } \|_1 \\
& \le  &\| \beta_{S^* } - \beta^* \|_1 - \| \beta_{-S^*} \|_1 \\
& = & ( \| b_{S^* } \|_1 - 
\| b_{-S} \|_1 ) \sigma^* / d .
\end{eqnarray*}
Therefore we obtain
$${\rm e}^{ P_n g_{0,1}} 
\biggl \{ P_n (g_{\hat b , \hat d } -  g_{0,1}) 
 + 1 - \hat d  \biggr \}   \le   \lambda ( \| \hat b_{S^*} \|_1 - \| \hat b_{-{S^* }}  \|_1) .$$
  But then also, using the convexity of $(b,d) \mapsto g_{b,d}$ and that $\hat b_t = t \hat b $ and 
 $1- \hat d_t = t (1- \hat d ) $, 
 $${\rm e}^{ P_n g_{0,1}} 
\biggl \{ P_n (g_{\hat b_t , \hat d_t } -  g_{0,1}) 
 + 1 - \hat d_t  \biggr \}   $$
 \begin{eqnarray*}
 & \le &{\rm e}^{ P_n g_{0,1}} \biggl \{ t P_n g_{\hat b , \hat d } + (1- t ) P_n   g_{0,1}- P_n g_{0,1} 
 + t (1- \hat d)  \biggr \}  \\
 & \le&  t \lambda ( \| \hat b_{S^*} \|_1 - \| b_{-S^* }\|_1 ) = \lambda ( \| \hat b_{t, S^* } \|_1 -
 \| \hat b_{t, -S^* } \|_1 ). 
 \end{eqnarray*} 
 \hfill$\sqcup \mkern -12mu \sqcap$
 
 \subsection{Excess risk}\label{excess-risk.section}
 
 Starting from the basic inequality for an empirical risk minimizer, a typical next step is
 to add and subtract the theoretical risk. In our case the theoretical risk
 is
 $$ R ( \beta , \sigma) = \EE R_n ( \beta , \sigma) , \
 (\beta , \sigma) \in \R^p \times \R_+ . $$ The true parameter $(\beta^* , \sigma^*)$
 is a minimizer of $R(\beta , \sigma)$, $(\beta, \sigma) \in \R^p \times R_+$, so under regularity
 the first derivative $\dot R ( \beta^* , \sigma^*) $ is zero. The excess risk at $(\beta , \sigma)$ is
 defined as
 $$   R ( \beta , \sigma) - R( \beta^* , \sigma^*) $$
 which is thus non-negative. 
 
 We use the re-parametrization $(\beta, \sigma) \mapsto (b,d)$.
 Then the excess risk with this new parametrization is
 $${\cal E} (b,d) := R ( \beta , \sigma) - R( \beta^* , \sigma^*) =  
 {1 \over n} \sum_{i=1}^n H(x_i b,d)  - \log d  =: {\cal H} (c,d) ,$$
 where
 $$ H(c,d) := \EE l( d \xi -c )  , \ c \in \R , \ d >0. $$
Note that ${\cal E} (b,d)$ is minimized at $(b^* , d^* )=(0,1)$.
 Moreover, under favorable conditions, for some constant $\kappa>0$, for $(b,d) $ in an appropriate neighborhood of $(0,1)$,
 $${\cal E} (b,d)  \ge {\kappa \over 2}  ( b^T \hat \Sigma b  + | d-1 |^2 ) . $$

 Instead of the basic inequality, we take the twisted basic inequality (\ref{twisted.equation})  as a starting point.
 Adding and subtracting the theoretical counterparts in
(\ref{twisted.equation})  gives
\begin{equation}\label{twisted2.equation}
 {\rm e}^{P_n g_{0,1} } \biggl \{ P ( g_{\hat b_t , \hat d_t} - g_{0,1} ) - \hat d_t +1 \biggr \}  $$ $$\le 
  - {\rm e}^{P_n g_{0,1} } \biggl \{(P_n - P)(g_{\hat b_t , \hat d_t} - g_{0,1} ) \biggr \}  + 
\lambda ( \| b_{t, S^*} \|_1   -\| b_{t, -S^*} \|_1 ) . 
\end{equation}
On the left-hand side of equation (\ref{twisted2.equation}) we now have what one might call the
``twisted excess risk" ${\cal E}_0 (b , d) $ at $(\hat b_t , \hat d_t )$, where
$$ {\cal E}_0 (b,d) : = {1 \over n} \sum_{i=1}^n H(x_i b,d) - d +1 .$$
 Instead lower-bounding the excess risk ${\cal E}(b,d)$ we now have to lower-bound ${\cal E}_0 (b,d)$.
 
 \begin{lemma} \label{H.lemma} Assume Condition \ref{H.condition}. 
   Then for
$ |d-1| + \max_{1 \le i \le n } | x_i b | $ small enough (depending only on $f$), it holds that
 $$ {\cal E}_0 ( b,d) \ge { \kappa_0^2 \over 4}  \biggl ( b^T \hat \Sigma b + (d-1)^2\biggr ) ,  $$
 where $\kappa_0^2 >0$ is the smallest eigenvalue of $\ddot H(0,1)$. 
 \end{lemma}
 
 {\bf Proof of Lemma \ref{H.lemma}.}
 Let 
 ${\cal H}(c,d) := H(c,d) - \log d $ and ${\cal H}_0 (c,d):= H(c,d) - d +1 $.  Then, with $c_i = x_i b$, $i=1 , \ldots , n$,  
 $$ {\cal E} (b,d) = {1 \over n } \sum_{i=1}^n {\cal H}( c_i, d) , \ 
 {\cal E}_0 (b,d) = {1 \over n } \sum_{i=1}^n {\cal H}_0( c_i, d) . $$ 
 Since ${\cal H} (c,d)$ is minimized at $(c,d) = (0,1)$ we see that 
 $\dot {\cal H}(0,1) =0 $. 
    But at $(c,d)=(0,1)$, 
 $  \dot {\cal H}_0 (0,1)=\dot {\cal H}(0,1)  $. 
 In other words ${\cal H}$ and ${\cal H}_0$ share the same stationary point $(c,d)= (0,1)$.
 Note moreover that $\ddot {\cal H}_0 (c,d) = \ddot H(c,d)$ for all $(c,d) \in \R \times R_+$. 
Condition \ref{H.condition} says that 
$\ddot H (c,d)  $ is continuous near $(0,1)$: for all $\epsilon >0$ there is a $\delta >0$ such that
$$ \|  \ddot H (c,d) - \ddot H(0,1) \|_{\infty}  \le \epsilon 
, $$ when $ | c | + | d-1|   \le \delta  $.
Now for $\max_{1 \le i \le 1 } | c_i |  + |d-1|   \le \delta $ 
we have, for all $i \in \{1 , \ldots , n \}$, for an intermediate point
$ \bar c_i = t c_i$ and $\bar d = t d + (1-t)$, $0 \le t \le 1$, of $(c_i,d)$ and $(0,1)$, 
$$ \| \ddot H(\bar c_i ,\bar d) - \ddot H (0,1)\|_{\infty} \le \epsilon  
, $$
so that by a two-term Taylor expansion with first term vanishing,  for all $i \in \{1 , \ldots , n \}$, 
\begin{eqnarray*}
 {\cal H}_0(c_i , d)  
& \ge &  { \kappa_0^2 \over 4}   \biggl ( c_i^2 + (d-1)^2\biggr  ),
\end{eqnarray*} 
if we take $ \epsilon = \kappa_0^2/8 $.
But then
\begin{eqnarray*}
  {\cal E} (b,d)  &=& {1 \over n} \sum_{i=1}^n {\cal H}_0 (c_i , d) 
 \ge  { \kappa_0^2 \over 4} \biggl  ( b^T \hat \Sigma b + (d-1)^2 \biggr ) . 
\end{eqnarray*}
\hfill $\sqcup \mkern -12mu \sqcap$

Instead of assuming Conditions \ref{G.condition} and \ref{H.condition}  we will now give sufficient conditions
to prove these.

\begin{lemma} \label{H2.lemma} Suppose $\ddot l (y)$ exists for all $y \in \R$, that $\ddot l (y) >0$ for all $y$, and that $\ddot l (\cdot)$ is  continuous. Assume moreover that for $M$ small enough,
$$ \ddot l (\xi + y ) \le G( \xi) , \ \forall \ | y | \le M (K_{\rm x} + \| \xi | )  , $$
where $\EE G^2 (\xi) < \infty$ and $\EE G^2 (\xi) \xi^4 < \infty$. 
Then Conditions \ref{G.condition} and  \ref{H.condition} are met.
\end{lemma}

{\bf Proof of Lemma \ref{H2.lemma}.} 
Condition \ref{G.condition} follows from
$$ |\dot l ( \xi+ y) - \dot l ( \xi + \tilde y ) | = \ddot l ( \xi + t y + (1- t ) \tilde y  ) | y - \tilde y | , $$
where $0\le t \le 1$. When $\max\{ | y | , |\tilde y |\} \le M ( K_{\rm x} + | \xi | )$, this is also
true for the intermediate point $t y + (1-t) \tilde y $ so then $\ddot l ( \xi + t y + (1- t ) \tilde y  ) \le G(\xi) $.

Set $h_{c,d} ( \xi) = l( d \xi - c)$.
We have
$$ \ddot h_{c,d}  (\xi) = \begin{pmatrix} \ddot l ( d \xi - c) & -\ddot l ( d \xi - c) \xi \cr
-\ddot l ( d \xi - c) \xi & \ddot l ( d \xi - c) \xi^2  \end{pmatrix} $$
So
$$ \EE \ddot h_{c,d} (\xi) = \begin{pmatrix} \EE \ddot l ( d \xi - c) & -\EE \ddot l ( d \xi - c) \xi \cr
-\EE \ddot l ( d \xi - c) \xi & \EE \ddot l ( d \xi - c) \xi^2  \end{pmatrix} $$
and in particular
$$\EE  \ddot h_{0,1} (\xi) =\begin{pmatrix} \EE \ddot l (  \xi ) & \-EE \ddot l (  \xi ) \xi \cr
-\EE \ddot l (  \xi ) \xi & \EE \ddot l (  \xi ) \xi^2  \end{pmatrix} .$$
By dominated convergence, for $| c| + |d-1| \le M $,
$$ \EE \ddot h_{c,d} (\xi) = \ddot H (c,d) $$
Write 
$$\gamma^2 := { \biggl (\EE \ddot l (  \xi  ) \xi \biggl )^2 \over  \EE \ddot l (  \xi ) \EE \ddot l ( \xi  ) \xi^2 } . $$
If $\gamma =\pm 1 $ we must have
$$ \sqrt {\ddot l ( \xi) } \xi  = C \sqrt {\ddot l ( \xi  )}  , \ {\rm almost \ surely }, $$
for some constant $C$. This is not the case because $\ddot l (\cdot) > 0 $ and $\xi$ is not constant. Therefore $| \gamma  | < 1 $.
It follows that 
$$\kappa_0^2 \ge ( 1- |\gamma||  ) \min \{ \EE \ddot l (  \xi ), \EE \ddot l (  \xi )\xi^2 \} . $$
Since $\ddot l $ is continuous,  we have for $| c| + |d-1| \rightarrow 0 $,
$$ \| \ddot h_{c,d} (z) - \ddot h_{0,1} (z)  \|_{\infty} \rightarrow 0 , \ \forall \ z \in \R.$$
By dominated convergence, then also 
$$ \| \ddot H (c,d) - \ddot H (0,1) \|_{\infty} = \| \EE  \ddot h_{c,d} (\xi) -  \EE \ddot h_{0,1} (\xi)  \|_{\infty} \rightarrow 0 . $$
\hfill $\sqcup \mkern -12mu \sqcap$

\subsection{Restricted eigenvalue}\label{restricted-eigenvalue.section}

Condition \ref{Sigma.condition} together with Condition \ref{s.condition} make it possible to lower
bound $b^T \hat \Sigma b / \| b \|_2^2$ for appropriate $b $.  The two conditions allow us to
conclude 
what  is known in the literature on the Lasso as the restricted eigenvalue condition (\cite{brt09}). Since we
need Condition \ref{s.condition} anyway in Theorem \ref{main.theorem}, we thought combining it with 
Condition \ref{Sigma.condition} is better than alternatively
imposing the restricted eigenvalue condition directly.

\begin{lemma}\label{restricted-eigenvalue.lemma}
Assume Condition \ref{Sigma.condition} and let
let $\eta^2 \gg K_{\rm x}^2 s_{\rm max} \sqrt{\log p / n}    / \Lambda_{\rm x}^2   $. Then 
for $n$ large enough, and for a vector $b \in \R^p$ satisfying
$\| b \|_1 \stackrel{\le}{\sim} \| b_{S^*} \|_1/ \eta $ we have
$$ b^T \hat \Sigma b\ge { \Lambda_{\rm x }^2 \over 2}  \| b \|_2^2 . $$
\end{lemma}

{\bf Proof of Lemma \ref{restricted-eigenvalue.lemma}.} 
This follows from
\begin{eqnarray*}
b^T \hat \Sigma b  &=& b^T \Sigma b + b^T ( \hat \Sigma- \Sigma ) b
 \ge
 \Lambda_{\rm x}^2 \| b \|_2^2 - \| \hat \Sigma- \Sigma \|_{\infty}
 \| b \|_1^2 \\
 & \stackrel{\ge} {\sim} & \Lambda_{\rm x}^2 \| b \|_2^2 - K_{\rm x}^2 \sqrt {\log p / n}   \| b \|_1^2 
 \stackrel{\ge} {\sim}  \Lambda_{\rm x}^2 \| b \|_2^2 -s^*  \sqrt {\log p / n }  K_{\rm x}^2  \| b_{S^*} \|_2^2/ \eta^2  
 \\
 & \ge &  \Lambda_{\rm x}^2 \| b \|_2^2 - {s_{\rm max} } \sqrt {\log p / n }  K_{\rm x}^2    \| b \|_2^2/ \eta^2 
  =  ( \Lambda_{\rm x}^2 - o(\Lambda_{\rm x}^2  ) )  \| b \|_2^2 
 \ge  { \Lambda_{\rm x}^2  \over 2} \| b \|_2^2 .
 \end{eqnarray*} 
 \hfill $\sqcup \mkern -12mu \sqcap$
 
 \subsection{Asymptotic continuity}\label{asymptotic-continuity.section}

Asymptotic continuity  plays an important role in empirical process theory.
It is about convergence to zero in probability, of the increments of an empirical process
$ \{ \sqrt n ( P_n - P) g : \ g \in {\cal G} \}$ indexed by a class of functions ${\cal G}$. In our case,
we look at the empirical process indexed by the gradients
$$ \biggl \{ \begin{pmatrix} \dot g_{b,d}^{\rm b} \cr \dot g_{b,d}^{\rm d}\cr \end{pmatrix} : \  (b,d) \in \Theta_M \biggr \} \subset \R^{p+1}  .$$ 
We normalize by $\sqrt {n / \log p}$ instead of $\sqrt n$ as this is the $\| \cdot \|_{\infty}$-rate of convergence
of $(P_n - P)$ at a fixed $(p+1)$-dimensional gradient vector. We need the asymptotic continuity in order to be able to show
that one can take the tuning parameter $\lambda$ close to $F^{-1} (1- \alpha) $. 

Let us start with upper bounding the random variable $\lambda^*$, to see that it is indeed
of order at most $ \sqrt {\log p / n}$ in probability. We phrase this in a more general context so that the result can also be applied elsewhere. The point of the next lemma is that a suitable sub-exponential tails condition is replaced
by an ``envelope condition" (see \cite{Duembgen:09} where this is further developed). 

For a $\{\varepsilon_i\}_{i=1}^{\infty} $ an i.i.d.\ Rademacher sequence 
(that is $\PP ( \varepsilon_i=1) = \PP (\varepsilon_i =-1) = 1/2$) independent of $\{ \xi_i\} _{i=1}^{\infty}$
we define
$$ P_n^{\varepsilon} g := {1 \over n } \sum_{i=1}^n \varepsilon_i g(x_i ,  \xi_i  ) . $$
The conditional expectation (probability) given $\vec{\xi}:= \{ \xi_i \}_{i=1}^{\infty} $ is written
as $\EE_{\vec{\xi} } $ ($\PP_{\vec{\xi} } $).

\begin{lemma}
  \label{lambda*.lemma} Let $\{ z_i \}_{i=1}^{\infty}$ be i.i.d.\ copies of a random variable ${\bf z} \in {\cal Z}$
  and let $\{ \psi_{i,j} : \ {\cal Z} \rightarrow \R  , \ 1 \le i \le n , \ 1 \le j \le p \} $ be a collection of real-valued
  functions on ${\cal Z}$, with $\EE \psi_{i,j}({\bf z} ) =0 $ for all $i$ and $j$,  and with envelope
  $$ \max_{1 \le j \le p } \max_{1 \le i \le n}  | \psi_{i,j}  (\cdot ) | \le \Psi (\cdot ) ,$$
  where 
  $$\EE \Psi^2  ({\bf z}) < \infty . $$
  Then
  $$ \max_{1 \le j \le p} {1 \over n} \sum_{i=1}^n \psi_{i,j} (z_i)  = {\mathcal O}_{\PP} ( \sqrt { \log p / n } ).  $$
  \end{lemma}
  
  {\bf Proof of Lemma \ref{lambda*.lemma}.} We consider the symmetrized version
  $$ {1 \over n} \sum_{i=1}^n \varepsilon_i \psi_{i,j} (z_i) $$
  with $\{ \varepsilon_i \} $ independent of $\{ z_i \}$. Then by Hoeffding's inequality 
  (\cite{Hoeffding:63}), for each $j$
  for all $t>0$, with probability at least $1- \exp[-t]$
  $$  \biggl | {1 \over n} \sum_{i=1}^n \varepsilon_i \psi_{i,j} (z_i) \biggr | \le \sqrt {8t} \sqrt{ {1 \over n} \sum_{i=1}^n \psi_{i,j}^2 }. $$
  Thus with probability at least $1- \exp[-t] $
  $$ \max_{1 \le j \le p}   \biggl | {1 \over n} \sum_{i=1}^n \varepsilon_i \psi_{i,j} (z_i) \biggr | \le \sqrt {8(t+ \log p)} \sqrt{ {1 \over n} \sum_{i=1}^n \Psi (z_i)^2 }. $$
  Therefore with probability at least $1- 1/p$,
  $$ \max_{1 \le j \le p}   \biggl | {1 \over n} \sum_{i=1}^n \varepsilon_i \psi_{i,j} (z_i) \biggr | \le \sqrt {16 \log p} \sqrt{ {1 \over n} \sum_{i=1}^n \Psi (z_i)^2 }. $$
  So with probability $1-o(1)$,
  $$ \max_{1 \le j \le p}   \biggl | {1 \over n} \sum_{i=1}^n \varepsilon_i \psi_{i,j} (z_i) \biggr | \le \sqrt {16 \log p)} \sqrt{2 \EE \Psi^2({\bf z} )} . $$
  But then, de-symmetrizing, with probability  $1- o(1)$,
  $$ \max_{1 \le j \le p}   \biggl | {1 \over n} \sum_{i=1}^n  \psi_{i,j} (z_i) \biggr | \le 4 \sqrt {16 \log p)} \sqrt{2 \EE \Psi^2({\bf z} )}. $$
  \hfill $\sqcup \mkern -12mu \sqcap$

\begin{theorem}\label{asymptotic-continuity.theorem} For a constant $C^{\rm b} $ depending only on $f$,
we have
$$ \sup_{(b,d) \in \Theta_M }
\biggl \| (P_n-P) (  \dot g_{b,d}^{\rm b}  -  \dot g_{0,1}^{\rm b}) \biggr \|_{\infty}  \le C^{\rm b}
K_{\rm x}^2 M \sqrt {\log p / n }  $$
with probability at least $1- 4/p- \alpha^{\rm b} $ where $\alpha^{\rm b} \rightarrow 0$ depends only on $f$.
Moreover,
$$ \sup_{(b,d) \in \Theta_M }
\biggl | (P_n-P) (  \dot g_{b,d}^{\rm d}  -  \dot g_{0,1}^{\rm d}) \biggr |  \le C^{\rm d}
K_{\rm x} M \sqrt {\log p / n }  ,$$
where the constant $C^{\rm d}$ and the sequence $\alpha^{\rm d} \rightarrow 0 $ depend only on $f$. 
\end{theorem}

{\bf Proof of Theorem \ref{asymptotic-continuity.theorem}.} 
Let
$$ \lambda_1^{\rm b} := { \EE_{\vec{\xi} } \biggl \| {1 \over n} \sum_{i=1}^n \varepsilon_i G( \xi_i) x_i \biggr \|_{\infty} 
\over K_{\rm x} \sqrt { {1 \over n} \sum_{i=1}^n G (\xi_i^2 ) } }, \ 
 \lambda_2^{\rm b} := { \EE_{\vec{\xi}} \biggl | {1 \over n} \sum_{i=1}^n \varepsilon_i G(\xi_i) \xi_i \biggr | \over
\sqrt { {1 \over n}  \sum_{i=1}^n G^2 (\xi_i) \xi_i^2 } } . $$
By Lemma \ref{lambda*.lemma}, we know that  $\lambda_1^{\rm b} \asymp  \sqrt {\log p / n } $. 
Moreover, $\lambda_2^{\rm b} \asymp  1/ \sqrt n  $. Invoking the contraction theorem 
(\cite{Ledoux:91}), we obtain for
$j\in \{1 , \ldots , p \}$
$$ \EE_{\vec{\xi} } \sup_{(b,d) \in \Theta_M }
\biggl | P_n^{\varepsilon} ( ( \dot g_{b,d}^{\rm b} )_j - ( \dot g_{0,1}^{\rm b})_j) \biggr | $$
$$ \le 2 K_{\rm x} M \biggl ( K_{\rm x} \lambda_1^{\rm b} \sqrt {{1 \over n} \sum_{i=1}^n G^2 ( \xi_i)} +
\lambda_2^{\rm b} \sqrt {{1 \over n} \sum_{i=1}^n G^2 ( \xi_i) \xi_i^2} \biggr ) . $$
We used here Condition \ref{G.condition} and the fact that 
$$ | ( \dot g_{b,d}^{\rm b}  ( x_i , \xi_i))_j - ( \dot g_{0,1}^{\rm b}(x_i , \xi_i))_j  | \le K_{\rm x}
| \dot l ( d \xi_i - b x_i) - \dot l (\xi_i) | . $$
Continuing with the latter we see that by again inserting Condition \ref{G.condition},
$$| ( \dot g_{b,d}^{\rm b}  ( x_i , \xi_i))_j - ( \dot g_{0,1}^{\rm b}(x_i , \xi_i))_j  | \le K_{\rm x} M 
G(\xi_i) (  K_{\rm x} + | \xi_i|   ) . $$
By Massart's concentration inequality (\cite{Massart:00a}), for all $j \in \{1 , \ldots , p \} $ and for all $t >0$, with $\PP_{\vec{\xi}}$ probability at least
$1- \exp[-t]$
$$ \sup_{(b,d) \in \Theta_M }
\biggl | P_n^{\varepsilon} ( ( \dot g_{b,d}^{\rm b} )_j - ( \dot g_{0,1}^{\rm b})_j) \biggr |  $$
\begin{eqnarray*}
& \le & 2 K_{\rm x} M \biggl ( K_{\rm x} \lambda_1^{\rm b} \sqrt {{1 \over n} \sum_{i=1}^n G^2 ( \xi_i)} +
\lambda_2^{\rm b} \sqrt {{1 \over n} \sum_{i=1}^n G^2 ( \xi_i) \xi_i^2} \biggr ) \\
& +&  K_{\rm x} M \biggl (  K_{\rm x} \sqrt { {1 \over n} \sum_{i=1}^n G^2 (\xi_i) } +
\sqrt{ {1 \over n} \sum_{i=1}^n G^2 (\xi_i) \xi_i^2 } \biggr ) \sqrt { 8t \over n} .
\end{eqnarray*}
Therefore, with $\PP_{\vec{\xi}}$ probability at least
$1- \exp[-t]$
$$ \sup_{(b,d) \in \Theta_M }
\biggl \| P_n^{\varepsilon} (  \dot g_{b,d}^{\rm b}  -  \dot g_{0,1}^{\rm b}) \biggr \|_{\infty}   $$
\begin{eqnarray*}
& \le & 2 K_{\rm x} M \biggl ( K_{\rm x} \lambda_1^{\rm b} \sqrt {{1 \over n} \sum_{i=1}^n G^2 ( \xi_i)} +
\lambda_2^{\rm b} \sqrt {{1 \over n} \sum_{i=1}^n G^2 ( \xi_i) \xi_i^2} \biggr ) \\
& +&  K_{\rm x} M \biggl (  K_{\rm x} \sqrt { {1 \over n} \sum_{i=1}^n G^2 (\xi_i) } +
\sqrt{ {1 \over n} \sum_{i=1}^n G^2 (\xi_i) \xi_i^2 } \biggr ) \sqrt { 8(t+\log p)  \over n} .
\end{eqnarray*}
The inequality also holds with $\PP$-probability at least $1- \exp[-t]$. We take $t=\log p $ and
let
$$ \alpha^{\rm b} / 4 := \PP \biggl (  \biggl \{ {1 \over n} \sum_{i=1}^n G^2 (\xi_i) > 2 \EE G^2 ( \xi) \biggr \} \cup 
\biggl \{ {1 \over n} \sum_{i=1}^n G^2 (\xi_i) \xi_i^2 > 2 \EE G^2 (\xi ) \xi^2 \biggr \} \biggr ) .$$
Note that $\alpha^{\rm b}$ depends only on $f$ and, by Condition \ref{G.condition}, $\alpha^{\rm b} \rightarrow 0 $.
Then with $\PP$-probability at least $1- 1/p -\alpha^{\rm b}/4$
$$ \sup_{(b,d) \in \Theta_M }
\biggl \| P_n^{\varepsilon} (  \dot g_{b,d}^{\rm b}  -  \dot g_{0,1}^{\rm b}) \biggr \|_{\infty}   $$
\begin{eqnarray*}
& \le & 2 K_{\rm x} M \biggl ( K_{\rm x} \lambda_1^{\rm b} \sqrt {2 \EE G^2 (\xi) } +
\lambda_2^{\rm b} \sqrt {2 \EE G^2 (\xi) \xi^2 } \biggr ) \\
& +&  K_{\rm x} M \biggl (  K_{\rm x} \sqrt {2 \EE G^2 ( \xi) } +
\sqrt{ 2 \EE G^2 (\xi) \xi^2 }  \biggr ) \sqrt { 16 \log p  \over n} \\
&=: &{C^{\rm b} \over 4} K_{\rm x}^2 M \sqrt {\log p \over n}  ,
\end{eqnarray*}
(say) where the constant $C^{\rm b}$ only depends on $f$.
De-symmetrizing gives that
$$ \sup_{(b,d) \in \Theta_M }
\biggl \| (P_n-P) (  \dot g_{b,d}^{\rm b}  -  \dot g_{0,1}^{\rm b}) \biggr \|_{\infty}  \le C^{\rm b}
K_{\rm x}^2 M \sqrt {\log p / n }  $$
with probability at least $1- 4/p - \alpha^{\rm b} $. \\
For the partial derivative with respect to $d$ we can use the same arguments. 
Define
$$ \lambda_1^{\rm d} := {\biggl \| {1 \over n} \sum_{i=1}^n \varepsilon_i G(\xi_i) \xi_i x_i \biggr \|_{\infty} 
\over K_{\rm x} \sqrt { {1 \over n } \sum_{i=1}^n G^2 (\xi_i) \xi_i^2 } } , \ 
 \lambda_2^{\rm d} := {\biggl | {1 \over n} \sum_{i=1}^n \varepsilon_i G(\xi_i) \xi_i^2 \biggr | \over 
 \sqrt { {1 \over n } \sum_{i=1}^n G^2 (\xi_i) \xi_i^4 } } , $$
 and
 $$ \alpha^{\rm d} / 4 := 
 \PP \biggl (  \biggl \{ {1 \over n} \sum_{i=1}^n G^2 (\xi_i) \xi_i^2 > 2 \EE G^2 ( \xi) \xi^2 \biggr \} \cup 
\biggl \{ {1 \over n} \sum_{i=1}^n G^2 (\xi_i) \xi_i^4 > 2 \EE G^2 (\xi ) \xi^4 \biggr \} \biggr ) .$$
We get
with probability at least $1- 4/p - \alpha^{\rm d} $
$$ \sup_{( b,d) \in \Theta_M } \biggl | (P_n - P) ( \dot g_{b,d}^{\rm d} - \dot g_{0,1}^{\rm d}) \biggr | $$
\begin{eqnarray*}
& \le & 8 M \biggl ( K_{\rm x} \lambda_1^{\rm d} \sqrt {2 \EE G^2 (\xi) \xi^2 } +
\lambda_2^{\rm d} \sqrt {2 \EE G^2 (\xi) \xi^4 } \biggr ) \\
& +& 4 M \biggl (  K_{\rm x} \sqrt {2 \EE G^2 ( \xi) xi_i^2 } +
\sqrt{ 2 \EE G^2 (\xi) \xi^4 }  \biggr ) \sqrt { 8 \log p  \over n} \\
&=: &C^{\rm d}  K_{\rm x} M \sqrt {\log p \over n}  .
\end{eqnarray*}
\hfill $\sqcup \mkern -12mu \sqcap$

 \subsection{Proof of Theorem \ref{main.theorem} using the preliminary results}\label{proof-main2.section}
 
 We present a more detailed version of Theorem \ref{main.theorem}.
 Define
  $$ \bar \kappa^2 :=  {\rm e}^{P g_{0,1} } {\kappa_0^2  \over 8}. $$
Note that $\bar \kappa^2$ is a constant depending on $f$ only, so it does not depend on $n$.
To facilitate checking the result, we however kept this constant explicitly  in the bounds.

  \begin{theorem} \label{main2.theorem}  Take $\lambda = {\mathcal O} (\sqrt {\log p / n})$, 
    $ \lambda (1- \eta)  \ge  F^{-1} (1- \alpha) $, 
  where $1- \eta \gg 1 $, 
  $$\eta^2 \gg { K_{\rm x}^2  \lambda s_{\rm max} \over \Lambda_{\rm x}^2 \bar \kappa^2}  +
  { K_{\rm x} \lambda \over \bar \kappa^2 } .$$  Let
  $$ {M \over 4} := {9 \over \eta}  {  8 \lambda  s^* \over  \Lambda_{\rm x}^2 \bar  \kappa^2 
   }  +{ 5 \over \eta} { 8 \lambda \over \bar \kappa^2 } . $$
  Then for $n$ large enough, we have with probability at least $1-\alpha + o(1)$.
  $$ \| \hat b \|_1 + | \hat d-1 | \le M , $$
  $$ \biggl ( \hat b^T \hat \Sigma  \hat  b \biggr )^{1/2}  \stackrel{\le} {\sim}
    {   \lambda  \sqrt { s^*}  \over  \Lambda_{\rm x} \bar  \kappa^2} + { \lambda  \over \bar \kappa^2 }  , $$
    $$ \| \hat b \|_2 \le {  \lambda  \sqrt { s^*}  \over  \Lambda_{\rm x}^2  \bar \kappa} +
    { \lambda \over \Lambda_{\rm x} \bar \kappa^2} 
    , $$
    and
    $$ | \hat d- 1| \stackrel{\le}{\sim}{\lambda  \sqrt {s^*} \over \Lambda_{\rm x} \bar \kappa^2 } +  {\lambda\over \bar \kappa^2 } . $$
   \end{theorem}

  {\bf Proof of Theorem \ref{main2.theorem}.}

  Recall that $F^{-1} (1- \alpha)$ is the $(1-\alpha)$-quantile of
  $$\lambda^*= {\rm e}^{P_n g_{0,1} } \| \dot g_{0,1}^{\rm b} \|_{\infty}.$$
  Thus, with probability $1- \alpha$, $\lambda^* \le F^{-1} ( 1-\alpha)  $.
  Since $\EE (\dot l (\xi ) \xi )^2 < \infty$ by Condition \ref{f.condition}, we know that for all $u>0$,
  $$ \PP \biggl ( |(P_n - P) \dot g_{0,1}^{\rm d} | \ge  \sqrt {u\log p / n } \biggr ) \le 
  {\EE (\dot l (\xi ) \xi )^2  \over u \log p  } \rightarrow 0 , \ p \rightarrow \infty 
    . $$
  By Condition \ref{alpha.condition}, we conclude that  with probability tending to 1, it holds that
  $$ {\rm e}^{P_n g_{0,1} }  |(P_n - P) \dot g_{0,1}^{\rm d} | \le  F^{-1} ( 1-\alpha) . $$
  By Theorem \ref{asymptotic-continuity.theorem}, with probability at least $1- \alpha^{\rm b} - \alpha^{\rm d}- 8/p$,
\begin{equation} \label{Mb.equation}
  \sup_{(b,d) \in \Theta_M } \biggl \| (P_n-P) (  \dot g_{b,d}^{\rm b}  -  \dot g_{0,1}^{\rm b}) \biggr \|_{\infty}  \le K_{\rm x}^2 \bar M  F^{-1} (1-\alpha) ,
 \end{equation}
and
\begin{equation}\label{Md.equation}
 \sup_{(b,d) \in \Theta_M } \biggl | (P_n-P) (  \dot g_{b,d}^{\rm d}  -  \dot g_{0,1}^{\rm d}) \biggr |  \le 
K_{\rm x} \bar M F^{-1} (1- \alpha)   , 
\end{equation}
where
$$ \bar M := { \max\{ C^{\rm b} , C^{\rm d} \} M \sqrt {\log p / n } \over F^{-1} (1- \alpha) } . $$
Note that since $M \rightarrow 0 $ and by Condition \ref{alpha.condition} , also $\bar M \asymp M  \rightarrow 0 $. 
In the rest of the proof, we assume we are on the set ${\cal A} $ where the above inequalities
  (\ref{Mb.equation}) and (\ref{Md.equation}) hold and where
  in addition
  $$ {\rm e}^{P_n g_{0,1} }  \max \biggl \{ \biggl \| (P_n - P) \dot g_{0,1}^{\rm b}\biggr \|_{\infty},  \biggl |(P_n - P) \dot g_{0,1}^{\rm d} \biggr | \biggr \}
  \le F^{-1} (1- \alpha) ,  \ P_n g_{0,1} \ge  { P g_{0,1} \over 2} . $$
  Note that $\PP ({\cal A} ) = 1- \alpha - o(1) $. 
  
   Define 
  $\hat b_t = t \hat b $ and $\hat d_t = t \hat d + (1- t) $
  where 
  $$t = {M \over M + \|\hat b\|_1 + |\hat d - 1|   } . $$
  Then
  \begin{eqnarray*}
  \| \hat b_t  \|_1 + | \hat d_t - 1 |   \le  M . 
  \end{eqnarray*} 
     By the mean-value theorem (in higher dimensions), and  interchanging
    integration and differentiation (which is allowed by dominated convergence in view of
    Condition \ref{G.condition}), for an intermediate point
 $  (\bar b , \bar d)$, we have 
 \begin{eqnarray*}
  & & 
\biggl |  \biggl ( P_n - P \biggr ) \biggl ( g_{\hat b_t , \hat d_t} - g_{0,1} \biggr ) \biggr |\\
& = &\biggl | (P_n-P) \begin{pmatrix}  \hat b_t   \cr
 \hat d_t-1  \cr \end{pmatrix}^T  \begin{pmatrix} \dot g_{\bar b , \bar d }^{\rm b} \cr
  \dot g_{\bar b , \bar d }^{\rm d} \cr \end{pmatrix}  \biggr | \\
  & \le& \biggl \| ( P_n - P  ) \dot g_{\bar b , \bar d  }^{\rm b} \biggr \|_{\infty} \| \hat b_t \|_1 +
  \biggl | ( P_n - P) \dot g_{\bar b , \bar d }^{\rm d} \biggr | | \hat d_t -1 |  .
  \end{eqnarray*}

  We apply the twisted basic inequality (\ref{twisted2.equation}), which yields that
  \begin{eqnarray*}
  & & {\rm e}^{P_n g_{0,1} } \biggl \{ P ( g_{\hat b_t , \hat d_t} - g_{0,1} ) - \hat d_t +1 \biggr \}  \\
  &\le &
  - {\rm e}^{P_n g_{0,1} } \biggl \{(P_n - P)(g_{\hat b_t , \hat d_t} - g_{0,1} ) \biggr \} +
\lambda ( \| \hat b_{t, S^*} \|_1   -\| \hat b_{t, -S^*} \|_1 ) \\ 
&\stackrel{{\rm on} \ {\cal A} } { \le}   &
F^{-1} (1-\alpha) \biggl \{ ( 1+ K_{\rm x}^2 \bar M ) \| \hat b_t \|_1 +
 (1+  K_{\rm x} \bar M) | \hat d_t -1 |  \biggr \}  \\ 
 &+ & 
\lambda ( \| \hat b_{t, S^*} \|_1   -\| \hat b_{t, -S^*} \|_1 ) \\
& \stackrel{ F^{-1} (1- \alpha ) \le \lambda (1- \eta)  } {\le} & 
\lambda (2-\eta+ K_{\rm x}^2 \bar M ) \| \hat b_{t, S^* }  \|_1 - \lambda (\eta - K_{\rm x}^2 \bar M ) 
\| \hat b_{t, - S^*} \|_1 \\
&+& \lambda (1- \eta) (1+  K_{\rm x} \bar M) | \hat d_t -1 |  \\
& \le & 2 \lambda  \| \hat b_{t, S^* }  \|_1 - { \eta \over 2}  \lambda  
\| \hat b_{t, - S^*} \|_1 + 2 \lambda  | \hat d_t -1 | 
\end{eqnarray*}   
where in the last step, we  invoke that for $n$ large enough, $\eta \ge 2 K_{\rm x}^2 \bar M $ and $K_{\rm x} \bar M \le 1$. 
Furthermore, by Lemma \ref{H.lemma}, for $n $ large enough
     $$ P(g_{\hat b_t , \hat d_t } -
     g_{0,1} ) +1 - \hat d_t   \ge  {\kappa_0^2  \over 4}  \biggl \{ \hat b_t^T \hat \Sigma 
   \hat b_t + |\hat d_t - 1 |^2  \biggr \} , $$
   since $\| \hat b_t \|_1 + | \hat d_t -1 | \le M = o(1)$ and so also $\max_{1 \le i \le n } | x_i \hat b_t | \le
   K_{\rm x} M = o(1) $.

   Thus we get
   \begin{equation}\label{cases.equation}
   {\rm e}^{P g_{0,1} } {\kappa_0^2  \over 4}  \biggl \{ \hat b_t^T \hat \Sigma 
   \hat b_t + |\hat d_t - 1 |^2  \biggr \}  + {\eta \over 2 } \lambda  \| \hat b_{t, - S^*} \|_1  
   \end{equation} 
   $$ \le \underbrace{2 \lambda  \| \hat b_{t, S^* }   \|_1}_{:= (i)} +  
   \underbrace{2\lambda  | \hat d_t -1 |}_{:=(ii) } .   $$
   Recall for the next arguments that
   $ \bar \kappa^2 := { {\rm e}^{P g_{0,1} } \kappa_0^2  /  8}   $.

     We now consider two cases:\\
     {\it Case 1.} $(i) \le (ii)$, \\
      {\it Case 2.} $  (i) \ge (ii) $.\\

        \underline {In Case 1} we find
   $$\bar \kappa^2 | \tilde d-1 |^2   \le  4 \lambda \ | \hat d_t -1 | ,$$
   or
   $$| \hat d_t- 1 |  \le  {4 \lambda \over  \bar \kappa^2  } \le { M\over 4} .$$
      Further
      \begin{eqnarray*}
    {\eta \over 2 } \| \hat b_t   \|_1   & = & {\eta \over 2} \| \hat b_{t, -S^* } \|_1  + {\eta \over 2} \| \hat b_{t, S^* } \|_1 \\
    & \le&  4 | \hat d_t -1| + { \eta } | \hat d_t -1 | 
     \stackrel {\eta \le 1}{\le}  {5} | \hat d_t-1| 
    \le { 5}  { 4 \lambda \over \bar \kappa^2} .  
   \end{eqnarray*} 
   This gives
   $$ \| \hat b_t \|_1 \le {5 \over \eta}  { 4 \lambda \over \bar \kappa^2 } \le { M \over 4} . $$
   So
   $$ \| \hat   b_t  \|_1 + | \hat d_t - 1 |   \le M/ 2 . $$
   Moreover, we get in Case 1,
   $$  \hat b_t^T \hat \Sigma \hat  b_t  \le 
    4  \lambda |\hat d_t -1|  \le  { (4 \lambda)^2 \over \bar \kappa^2 } . $$
       But then
    \begin{eqnarray*}
    \hat b_t^T  \hat \Sigma \hat b_t &\ge& { \Lambda_{\rm x}^2  } 
    \| \hat b_t^T \|_2 - \| \hat \Sigma - \Sigma \|_{\infty} \| \hat b_t \|_1^2 \\
    &=& 
   \Lambda_{\rm x}^2  \| \hat b_t^T \|_2^2 - {\mathcal O} ({K_{\rm x}^2 \lambda^3  \over \eta^2 \bar \kappa^4 } )  
   = \Lambda_{\rm x}^2  \| \hat b_t^T \|_2^2 -  {  \Lambda_{\rm x}^2 \lambda^2 \over \bar \kappa^2 } 
   o(1) , 
   \end{eqnarray*}
   since we assume $\eta^2 \gg K_{\rm x}^2 \lambda / (\Lambda_{\rm x}^2\bar \kappa^2 )   $ and
   $\| \Sigma - \hat \Sigma \|_{\infty} \stackrel {\le } {\sim} K_{\rm x}^2\lambda $. 
    So  
    \begin{eqnarray*}
    \| \hat b_t \|_2^2 &\le& {1 \over \Lambda_{\rm x}^2}  \biggl ( { \hat b_t^T  \hat \Sigma \hat b_t  } + 
    {\lambda^2 \over \bar \kappa^2 } o (1)  ) \biggr ) 
     \stackrel{\le}{\sim} 
     {  \lambda^2 \over \Lambda_{\rm x}^2   \bar \kappa^2 } . 
 \end{eqnarray*}

   \underline {In Case 2} we see that
   \begin{eqnarray*}
   { \eta \over 2}  \| \hat   b_t  \|_1 &= & {\eta \over 2} \| \hat b_{t, -S^* } \|_1  + {\eta \over 2} \| \hat b_{t, S^* } \|_1  
    \stackrel{\eta \le 1 } {\le}  { 9 \over 2}  \| \hat b_{t, S^*} \|_1  ,
   \end{eqnarray*}
   or
   $$ \| \hat b_t \|_1 \le {9 \over \eta} \| \hat b_{t, S^*} \|_1 . $$
   This implies by Lemma \ref{restricted-eigenvalue.lemma}, for $n$ large enough
   $$  \hat b_t^T \hat \Sigma \hat b_t \ge { \Lambda_{\rm x}^2 \over 2} \| \hat b_t\|_2^2  . $$
   We arrive at
   \begin{eqnarray*}
     {\bar \kappa^2 } \hat b_t^T \hat \Sigma \hat  b_t 
   & \le & 4 \lambda    \| \hat  b_{t, S^*}   \|_1  
    \le  4 \lambda  \sqrt {s^*}  \| \hat  b_t  \|_2  
    \le  { 4 \over \Lambda_{\rm x} } \lambda \sqrt {2 s^*} 
       \biggl (  \hat b_t^T \hat \Sigma  \hat  b_t  \biggr )^{1/2}  .
       \end{eqnarray*}
    So
    $$ \biggl ( \hat b_t^T \hat \Sigma  \hat  b_t \biggr )^{1/2}  \le
    { 4  \lambda  \sqrt { 2s^*}  \over  \Lambda_{\rm x}  \bar \kappa^2} . $$
    Then also
    $$ \| \hat b_t \|_2 \le { \sqrt {2} \over \Lambda_{\rm x} } ( \hat b_t^T \hat \Sigma \hat b_t )^{1/2} 
    \le { 8 \lambda {\sqrt s_*} \over \Lambda_{\rm x}^2 \bar \kappa^2 } . $$
    But then
    \begin{eqnarray*}
     \| \hat b_{t}  \|_1 &\le &{9 \over \eta}  \| \hat b_{t, S^*} \|_1 \le { 9 \over \eta}  \sqrt {s^*} \| \hat b_t \|_2 
    \le  { 9 \over \eta }  {  8 \lambda  s^*\over \Lambda_{\rm x}^2 \bar \kappa^2 } \le { M \over 4} .
    \end{eqnarray*}
    Moreover, as we are in Case 2, also 
   $$ | \hat d_t - 1 | \le  \| \hat b_{t, S^* } \|_1 \le 
  \sqrt {s^*} \| \hat b_t \|_2 \le { 8 \lambda s_*   \over \Lambda_{\rm x}^2 \bar \kappa^2 } \le { M \over 4}  .$$
   Therefore, also in Case 2, 
 $$  \| \hat  b_t \|_1 + |\hat d_t -1  | \le M/2 . $$
 In fact, from (\ref{cases.equation}),
  $$ \bar \kappa^2 | \hat d_t-1 |^2 \le 4 \lambda \| \hat b_{t, S^*} \|_1 \le 4 \lambda \sqrt {s^*}
  \| \hat b_{t, S^* } \|_2 \le 4 \lambda \sqrt {s^*} { 8 \lambda \sqrt {s^*} \over \Lambda_{\rm x}^2 \bar \kappa^2 } . $$
 Because in both Case (i) and Case (ii), the bound $  \| \hat  b_t \|_1 + |\hat d_t -1  | \le M/2 $ is true, we 
 have now shown that
 $$ \ \| \hat  b\|_1 + |\hat  d -1  |  \le M . $$
 We can redo the proof with $(\hat b_t , \hat d_t)$ replaced by
 $( \hat b , \hat d)$.  
\hfill $\sqcup \mkern -12mu \sqcap$

\begin{corollary} \label{main1.corollary} Since with probability at least $1- \alpha + o(1)$, 
$| \sigma^* / \hat \sigma  -1 | = | \hat d -1 | \ll 1$ we see that 
with probability at least $1- \alpha + o(1)$, also
$$ \biggl ( (\hat \beta - \beta^*)^T \hat \Sigma ( \hat \beta - \beta^* ) \biggr )^{1/2} / \sigma^*
 \stackrel{\le} {\sim}
    {   \lambda  \sqrt { s^*}  \over  \Lambda_{\rm x} \bar  \kappa^2} + { \lambda  \over \bar \kappa^2 }  , $$
    and
    $$ \| \hat \beta - \beta^*  \|_2 / \sigma^*  \stackrel{\le}{\sim} {  \lambda  \sqrt { s^*}  \over  \Lambda_{\rm x}^2  \bar \kappa} +
    { \lambda \over \Lambda_{\rm x} \bar \kappa^2} 
    .$$
\end{corollary}

\section{Further proofs} \label{further-proofs.section} 

\subsection{Proof of the results in Section \ref{phase+.section}}
We start with an expansion of $\exp[R_n ( \beta , \sigma)] /\sigma^*= \exp[P_n g_{b,d} - \log d ] $
for small values of $\|b\|_1 $ and $d-1$. This will be applied in Theorem \ref{variable-selection.theorem}
for variable selection, 
and in Lemmas \ref{phase-transition1.lemma} and \ref{phase-transition2.lemma} for the result on the detection edge. Recall we assumed without loss of generality that $\sigma^*=1$.

\begin{lemma}\label{local-expansion.lemma}
Suppose the conditions of  Theorem \ref{main.theorem} with $\Lambda_{\rm x} $ fixed, and in addition Condition \ref{H2.condition}. Let $\sum_{i=1}^n x_{i, -0} =0 $. 
Consider sequences
$M^* =  {\mathcal O} (\lambda s^* / \eta )$, where  $K_{\rm x}^2 M^* \le \eta^2 $,  and
$e^* = \mathcal O( \lambda \sqrt {s^*}) $. 
We have with probability tending to 1, uniformly for all $( b_0 , \tilde b_{-0} , d) $ and $( b_0 , b_{-0} , d )$ in the set
$$ \Theta_{\rm local} := \biggl \{ ( b_0 , b_{-0} , d ) : \ \| b_{-0} \|_1 \le M^* , |d-1 | + |b_0| + \sqrt {b_{-0}^T \hat \Sigma_{-0} b_{-0} } \le e^* \biggr \} , $$
the local expansion
\begin{eqnarray*}
& & \exp [ P_n g_{b_0 , \tilde b_{-0} , d} ] - \exp[ P_n g_{b_0 , b_{-0} , d} ] \\
& = &  (   1+ {\mathcal O} ( \lambda M^*)  )\exp[ P_n g_{0,0,1} ]
 P_n \dot g_{0,0,1}^{{\rm b}_{-0}^T} ( \tilde b_{0} - b_{-0) }   +{\mathcal O} ( K_{\rm x}^2 M^*\lambda  )  
 \| \tilde b_{-0} - 
b_{-0} \|_1   \\ & +& \ddot {\cal H}^{{\rm c}, {\rm c}} (0,0,1) \bar b_{-0}^T \hat \Sigma_{-0} (\tilde b_{-0} - b_{-0}  ),
\end{eqnarray*} 
where $\bar b_{-0}= t \tilde b_{-0} + (1-t) b_{-0}$ ($0 \le t \le 1$) is an appropriate  intermediate point of $\tilde b_{-0}$ and $b_{-0} $. 
\end{lemma}

{\bf Proof of Lemma \ref{local-expansion.lemma}.} 
By an application of Theorem \ref{asymptotic-continuity.theorem}, we see that with probability tending to
1, uniformly for all $(b,d)$ in the set
$ \Theta_{\rm local} $, the validity of the bounds
\begin{eqnarray*}
 & &(P_n - P) (g_{b, d} - g_{0,1} ) = 
\underbrace{\| (P_n - P ) \dot g_{0,1}^{\rm b} \|_{\infty} \| b \|_1}_{={\mathcal O} 
( \lambda M^* ) } \\
&+& \underbrace{(P_n -P)\dot g_{0,1}^{\rm d} (d-1) }_{=  {\mathcal O} (\lambda e^* )} + \underbrace{{\mathcal O} (\lambda M^{*2} ) }_{= {\mathcal O} (e^{*2}  )} .
 \end{eqnarray*}
 Furthermore,
 \begin{eqnarray*}
  & &P (g_{b, d} - g_{0,1} - \log d )   =  \begin{pmatrix} d-1 \cr b_0 \cr  \end {pmatrix}^T
 \biggl (  \ddot {\cal H} (0,0,1)  +o(1) \biggr ) 
 \begin{pmatrix} d-1 \cr b_0 \cr  \end {pmatrix}^T  \\
 &+&
 \biggl ( \ddot {\cal H}^{{\rm c}, {\rm c} } (0,1)+o(1) \biggr )  \biggl ( b_{-0}^T \hat \Sigma_{-0} b_{-0} \biggr ) 
  \\
 & =& {\mathcal O} (e^{*2} ) . 
 \end{eqnarray*} 
 Thus
 $$ P_n ( g_{b,d} - g_{0,1} - \log d)  = {\mathcal O} ( \lambda M^*  )  . $$
 It follows that
 $$ \exp [ P_n g_{b,d}  ] - \exp [ P_n g_{0,1} ] = \exp [ P_n g_{0,1} ]  ( 1+  
 {\mathcal O} ( \lambda M^*)  ) .$$
 For $(b,d)$ and $(\tilde b , \tilde d )$ in $\Theta_{\rm local}$,  with an intermediate point
 $(\bar b, \bar d):= t(\tilde b , \tilde d ) +(1-t) ( b,d) $, also
 $$ \exp [ P_n g_{\bar b,d}  ] - \exp [ P_n g_{0,1} ] = \exp [ P_n g_{0,1} ]  ( 1+  
 {\mathcal O} ( \lambda M^*)  ) .$$
 But then
 \begin{eqnarray*}
 &  &\exp [ P_n g_{\tilde b,\tilde d} - \log \tilde d] - \exp[P_n g_{b,d} - \log d ]\\
 &=&
 \exp [ P_n g_{\bar b, \bar d} ]  \biggl ( P_n \dot g_{\bar b , \bar d }^{{\rm b}^T } ( \tilde b - b )
  + 
 P_n \dot g_{\bar b , \bar d }^{\rm d}  ( \tilde d-d) - { \tilde d - d \over \bar d} \biggr ) \\
 &= & \exp [ P_n g_{0,1} ] (1+ {\mathcal O} ( \lambda M^* ) ) 
 \biggl (  P_n \dot g_{\bar b , \bar d }^{{\rm b}^T } ( \tilde b -  b ) +
 P_n \dot g_{\bar b , \bar d }^{\rm d}  ( \tilde d-d) - { \tilde d - d \over \bar d}  \biggr ) \\
 &=&  P_n \dot g_{\bar b , d }^{{\rm b_{-0} }^T } ( \tilde b_{-0} -  b_{-0}  ) ,
 \end{eqnarray*}
 where the last equality is true when  $\tilde b_0= b_0$ and $\tilde d = d $.
  By Theorem \ref{asymptotic-continuity.theorem}, with probability tending to 1,
  $$(P_n - P)\dot g_{b_0, \bar b_{-0} , d }^{{\rm b}_{-0} ^T } ( \tilde b_{-0} -  b_{-0} ) =
  (P_n -P) \dot g_{0,0,1}^{{\rm b}_{-0}^T}  ( \tilde b_{-0} - b_{-0} )+ {\mathcal O} ( K_{\rm x}^2 M^* \lambda ) \| 
  \tilde b_{-0} - b_{-0} \|_1 .$$
   Moreover, for a further intermediate point $\bar {\bar b} $, with 
   $\dot {\cal H}^{\rm c}$ the derivative of ${\cal H}$ with respect to $c$,  
   \begin{eqnarray*}
  & & P \dot g_{b_0, \bar b_{-0} , d }^{{\rm b}_{-0}  }=  
   {1 \over n} \sum_{i=1}^n \dot {\cal H}^{\rm c} ( b_0 , x_{i, -0} \bar b_{-0} , d) ) x_{i, -0}^T \\
  & = & \underbrace{\dot {\cal H}^{\rm c} (0,0,1) }_{=0} {1 \over n } \sum_{i=1}^n x_{i,-0}^T + 
  { 1 \over n} \sum_{i=1}^n \ddot {\cal H} ({\bar {\bar b}}_0 + x_{i, -0} {\bar {\bar b}}_{-0} , {\bar {\bar d} } )  x_{i, -0}^T \begin{pmatrix} d-1 \cr b_0 \cr \end{pmatrix}  \\ 
  & + & {1 \over n} \sum_{i=1}^n \dot {\cal H}^{{\rm c}, {\rm c}} ( {\bar {\bar b}}_0 + x_i {\bar {\bar b}}_{-0} , {\bar {\bar d} } ) x_{i, -0}^T x_{i, -0}  \bar b_{-0}  .
  \end{eqnarray*}  
  But
  $$  \| \ddot {\cal H} ({\bar {\bar b}}_0 , x_{i,-0} {\bar {\bar b}}_{-0} , {\bar {\bar d} } ) 
  - \ddot {\cal H} (0,0,1) \|_{\infty} \le L_{\rm H} ( | {\bar {\bar b}}_0 + x_{i, - 0} {\bar {\bar b}}_{-0}| +  | {\bar {\bar d} }-1| ) .$$ 
  So, invoking that $\sum_{i=1}^n x_{i,j}= 0$ for $j \ge 2$,
 $$ \biggl \|  { 1 \over n} \sum_{i=1}^n \ddot {\cal H} ({\bar {\bar b}}_0 + x_{i,-0} {\bar {\bar b}}_{-0} , {\bar {\bar d} } )  x_{i,-0}^T   \biggr \|_{\infty} 
 \le  L_{\rm H} \underbrace{ {1 \over n} \sum_{i=1}^n (| {\bar {\bar b}}_0 + x_{i , -0}{\bar {\bar b}}_{-0}| +  | {\bar {\bar d} }-1| ) K_{\rm x} }_{
 = {\mathcal O} ( K_{\rm x}  e^*)}   . $$
 Thus
 $$\biggl | { 1 \over n} \sum_{i=1}^n \ddot {\cal H} ({\bar {\bar b}}_0 + x_{i,-0} {\bar {\bar b}}_{-0} , {\bar {\bar d} } )  x_i^T \begin{pmatrix} d-1 \cr b_0 \cr \end{pmatrix} \biggr |   = {\mathcal O} 
 ( K_{\rm x} e^{*2}  ) . $$
 Moreover 
 $$ \biggl \| {1 \over n} \sum_{i=1}^n \ddot {\cal H}^{{\rm c}, {\rm c}} ( {\bar {\bar b}}_0 + x_{i-0} {\bar {\bar b}}_{-0} , {\bar {\bar d} } ) x_{i, -0}^T x_{i , -0}  \bar b_{-0}  - \ddot {\cal H}^{{\rm c}, {\rm c}} (0,0,1) \hat \Sigma_{-0} \bar b_{-0} \biggr \|_{\infty}  $$
 $$ \le  L_{\rm H}  {1 \over n} \sum_{i=1}^n( | {\bar {\bar b}}_0 + x_{i,-0} {\bar {\bar b}}_{-0}| +  | {\bar {\bar d} }-1| ) 
 |x_{i,-0}  \bar b_{-0} | K_{\rm x}   = {\mathcal O} ( K_{\rm x} e^{*2}  ) .  $$
 We conclude that
 \begin{eqnarray*}
 & &P_n \dot g_{b_0, \bar b_{-0} , d}^{{\rm b}_{-0}T } ( \tilde b_{-0} - b_{-0} )\\
 &  =&
 P_n \dot g_{0,0,1}^{{\rm b}_{-0}^T} (\tilde b_{-0} - b_{-0} ) + \underbrace{{\mathcal O} ( K_{\rm x}^2 M^*\lambda  ) + 
{\mathcal O} (K_{\rm x} e^{*2} )}_{={\mathcal O} ( K_{\rm x}^2 M^* \lambda  )}  \|  \tilde b_{-0} - 
b_{-0} \|_1   \\ 
& +&  \ddot {\cal H}^{{\rm c}, {\rm c}} (0,0,1) \bar b_{-0}^T \hat \Sigma_{-0} (\tilde b_{-0} - b_{-0}  )^T.
\end{eqnarray*} 
   \hfill $\sqcup \mkern -12mu \sqcap$
   
   {\bf Proof of Theorem \ref{variable-selection.theorem}.} 
 Take
 $$ \hat a_{S^*} := (X_{S^*}^T X_{S^*} )^{-1} X_{S^*}^T X_{-S^*} \hat b_{-S^*} , $$
 and
 $$ \hat \alpha_{-S} = \hat a_{S^*} / \hat d =
 (X_{S^*}^T X_{S^*} )^{-1} X_{S^*}^T X_{-S^*} \hat \beta_{-S^*} . $$

We apply Lemma \ref{local-expansion.lemma} with $$ \tilde b_{-0} :=  \begin{pmatrix}
\hat b_{S^*} \cr \hat b_{-S^*} \cr \end{pmatrix} $$ and
$$ b_{-0} := \begin{pmatrix}
\hat b_{S^* } + \hat a_{S^* }   \cr 0  \cr \end{pmatrix}  . $$
Then by the irrepresentable condition  $\| \hat a_{S^*} \|_1 \le \| \hat b_{-S^*} \|_1$.
 Moreover, since projecting a vector cannot increase its length, 
 $$ \| X_{S^*} \hat a_{S^*}\|_2  \le 
 \| X_{-S^*} \hat b_{-S^*} \|_2 .  $$
  Therefore, with probability at least $1- \alpha + o(1)$, with this choice of $\tilde b_{-0}$ and $b_{-0}$, we are in
  $\Theta_{\rm local}$.

Then for $\bar b_{-S^*}$ an intermediate point of $\hat b_{-S^*}$ and $0$
$$ \bar b_{-0}^T \hat \Sigma_{-0} (\tilde b_{-0} - b_{-0}  ) = \bar b_{-S^*}^T Z_{-S^*}^T Z_{-S^*} \hat b_{-S^*} \ge 0 . $$
Therefore, by Lemma \ref{local-expansion.lemma}, with probability at least $1- \alpha + o(1)$,
\begin{eqnarray*}
& & \exp [ P_n g_{b_0 , \hat b_{S^*}, \hat b_{-S^*}  , d} ] - \exp[ P_n g_{b_0 ,\hat b_{S^*} + \hat a_{S^*}  , 0, d} ] \\
& \ge &  - \biggl ((   1+ {\mathcal O} ( \lambda M^*)  )\lambda^* + {\mathcal O} ( K_{\rm x}^2 M^*\lambda)\biggr ) 
 (\| \hat a_{S^*} + \hat b_{-S^*} \|_1 ) \\
 & \ge & - \lambda (1- r_{n,1} \eta)(1- \eta)  (\| \hat a_{S^*} \|_1  + \| \hat b_{-S^*} \|_1 ) ,
  \end{eqnarray*} 
  where $r_{n,1} = o(1)$. 
  Here we used that $\eta^2 \gg K_{\rm x}^2 M^* $.  Next, we have, when $| \hat d-1| = \mathcal O( e^*)$
  $$ \lambda \| \hat \beta_{-S^*} \|_1 \ge \lambda (1- \eta r_{n,2}) ) \| b_{-S^*} \|_1 $$
  with $r_{n,2} = o(1)$, 
  and
  $$ \lambda ( \| \hat \beta_{S^*} \|_1 - \| \hat \beta_{S^*} + \hat \alpha_{S^*} \|_1)  \ge
   - \lambda \| \hat \alpha_{S^*} \|_1 \ge - \lambda (1+ r_{n,3}  \eta) \| \hat a_{S^*} \|_1  $$
   with $r_{n,3}= o(1)$.
   So
   $$ - \lambda (1- \eta) (1- r_{n,2} ) (\| \hat a_{S^*} \|_1  + \| \hat b_{-S^*} \|_1 ) + \lambda \| \hat \beta_{S^*} \|_1 +
   \lambda \| \hat \beta_{- S^* }\|_1 - \lambda \| \hat \beta_{S^*} + \hat \alpha_{S^*} \|_1 $$ $$
   \ge - \lambda (1- \eta - r_{n,1} \eta ) (\| \hat a_{S^*} \|_1  + \| \hat b_{-S^*} \|_1 )  + 
   \lambda (1- r_{n,2} ) \eta  ) \| b_{-S^*} \|_1 - \lambda (1+ r_{n,3}  \eta) \| \hat a_{S^*} \|_1$$
   $$ = \lambda (\eta  - (r_{n,1} + r_{n,2} ) \eta  \| \hat b_{-S^* } \|_1  + \lambda (2-\eta - (r_{n,1} + r_{n,3} ) \eta )   \| \hat a_{S^* } \|_1 > 0 $$
   for 
   $$\| \hat a_{S^*} \|_1 < { \eta - (r_{n,1} + r_{n,2} )\eta   \over 2 - \eta - (r_{n,1} + r_{n,3} ) \eta } \| \hat b_{-S^* } \|_1 . $$
   \hfill $\sqcup \mkern -12mu \sqcap$
   
   {\bf Proof of Lemma \ref{phase-transition1.lemma}.}
We may apply Lemma \ref{local-expansion.lemma} with $M^* \asymp e^* \asymp \lambda$, with
$\tilde b_{-0}= \hat \beta_{-0}$ and $b_{-0} = 0$. We then do not need Condition \ref{H2.condition} but
apply that since
$\exp[ P_n g_{ b_0, b_{-0} , d}  - \log d ] $ is convex in $b_{-0}$ ,
\begin{eqnarray*}
& & \exp[ P_n g_{ \hat b_0, b_{-0} , \hat d} - \log \hat d ]  - \exp[ P_n g_{\hat b_0 , 0 , \hat d} - \log \hat d   ] \\
& \ge &
\exp[P_n g_{\hat b_0 , 0 , \hat d }   - \log d ] P_n 
(\dot g_{\hat b_0 , 0, \hat d}^{{\rm b}_{-0}} )^T  \hat b_{-0} .
\end{eqnarray*}
Then we apply that
$$ P (\dot g_{\hat b_0, 0, \hat d}^{\rm b_{-0}} - \dot g_{0,0,1}^{\rm b_{-0} } )=0,$$
where we used that $\sum_{i=1}^n x_{i,j} = 0$ for $j \in \{ 2 , \ldots , p \} $.
\hfill $\sqcup \mkern -12mu \sqcap$

{\bf Proof of Lemma \ref{phase-transition2.lemma}.}

By Theorem \ref{asymptotic-continuity.theorem}, and with the notation used there,
with probability 
$1- o(1) $
$$ \sup_{(b,d) \in \Theta_M }
\biggl \| (P_n-P) (  \dot g_{b,d}^{\rm b}  -  \dot g_{0,1}^{\rm b}) \biggr \|_{\infty}  \le C^{\rm b}
K_{\rm x}^2 M \sqrt {\log p / n }  ,$$
and
$$ \sup_{(b,d) \in \Theta_M }
\biggl | (P_n-P) (  \dot g_{b,d}^{\rm d}  -  \dot g_{0,1}^{\rm d}) \biggr |  \le C^{\rm d}
K_{\rm x} M \sqrt {\log p / n }  .$$
We place ourselves on the set ${\cal A}$ where the above two inequalities hold, where
$P_n g_{0,0,1} \ge  P g_{0,0,1} /2 $, and where $\lambda^* \stackrel{\le}{\sim} \sqrt {\log p / n} $.
We now want to also assume that $| \hat d - 1| \stackrel{\le}{\sim} \sqrt{\log p / n}$ and
$|\hat b_0 | \stackrel {\le}{\sim} \sqrt {\log p / n } $. 
It is easy to see that this is the case when $\hat b_{-0} =0$. Since the
$\hat b_{-0} \not= 0 $ is what we aim at proving, we from now on assume
that indeed  $| \hat d - 1| + | \hat b_0| \stackrel{\le}{\sim} \sqrt{\log p / n}$.
 We add these events and the event $\lambda\le \lambda^* (1- \eta) $ to our set ${\cal A}$, where we
 take $M= {\mathcal O} (\sqrt {\log p / n})$. 
Recall that $\lambda_{-0}^* = \exp[P_n g_{0,0,1} ] \| P_n \dot g_{0,0,1}^{{\rm b}_{-0} } \|_{\infty}$.

Let $\lambda_j^* := 
 := \exp[P_n g_{0,0,1 } ] ( P_n \dot g_{0,0,1}^{{\rm b}_{-0} })_j $, $j=2, \ldots , p$ and
$|\lambda_{\bf j}^* |= 
\max_{2 \le j \le p} | \lambda_j^* | $. Define
$$ \bar \kappa^2 := {\cal H}^{{\rm c} , {\rm c} } (0,1) \exp[P_n g_{0,0,1} ]  . $$Take
$(\tilde b_{-0} )_j=0 $ for $1< j \not= {\bf j} $ and 
$$(\tilde b_{-0})_{\bf j} = \begin{cases} - {  \lambda_{\bf j}^*- \lambda   \over 2 \bar \kappa^2 \| x_{\bf j} \|_2^2 /n }
 &
\lambda_{\bf j}^*  > 0 \cr  - {  \lambda_{\bf j}^*+ \lambda   \over \bar \kappa^2 \| x_{\bf j} \|_2^2 /n }
 & 
 \lambda_{\bf j}^* < 0 \cr \end{cases}  . $$
 Then
 $$ \| \tilde b_{-0} \|_1=  { |  \lambda_{\bf j}^* - \lambda  | \over 
 2 \bar \kappa^2 \| x_{\bf j}  \|_2^2 /n } = {\mathcal O} (   \sqrt {\log p / n } ) .  $$
  
We have for an intermediate point $\bar b_{-0} = t \tilde b_{-0} $, $0 \le t \le 1$, 
\begin{eqnarray*}
& & \exp[ P_n g_{\hat b_0 , \tilde b_{-0} , \hat d } - \log \hat d ] -\exp[ P_n g_{\hat b_0 , 0 , \hat d } - \log \hat d ] \\
&=&
\exp[ P_n g_{\hat b_0 , \bar b_{-0} , \hat d } - \log \hat d ]  P_n 
(\dot g_{\bar b_0 , \bar b_{-0}, \bar d }^{{\rm b}_{-0}} )_{\bf j}  (\tilde b_{-0} )_{\bf j} . 
\end{eqnarray*} 
Using the same arguments as in Lemma \ref{local-expansion.lemma} we obtain
$$ \exp [ P_n  g_{\hat b_0 , \bar b_{-0} , \hat d } - \log \hat d ]  -
\exp[P_n g_{0,0,1} ] = \exp[P_n g_{0,0,1} ] (1+ {\mathcal O} (\log p / n ) ) $$
We further have
$$P_n 
(\dot g_{\bar b_0 , \bar b_{-0}, \bar d }^{{\rm b}_{-0}} )_{\bf j}    
 = P_n  ((\dot g_{\bar b_0 , \bar b_{-0}, \bar d }^{{\rm b}_{-0}} )_{\bf j}    - 
(\dot g_{0,0,1}^{{\rm b}_{-0}} )_{\bf j} ) ) + P_n (\dot g_{0,0,1}^{{\rm b}_{-0}} )_{\bf j} ) )  ,$$
and
$$ \biggl \| (P_n - P)  
((\dot g_{\hat b_0 , \bar b_{-0}, \hat d}^{{\rm b}_{-0}} )_{\bf j} - (\dot g_{0,0,1}^{{\rm b}_{-0}} )_{\bf j} )  \biggr \|_{\infty}  =
{ \mathcal O} (K_{\rm x}^2 \log p / n ) . $$
Furthermore by Condition \ref{H2.condition},
$$P (\dot g_{\bar b_0 , \bar b_{-0}, \bar d }^{{\rm b}_{-0}} )_{\bf j}  =
\underbrace{{\mathcal O} (\sqrt {\log p / n } )^T\begin{pmatrix} \bar d-1 \cr \bar b_0 \cr \end{pmatrix} }_{={\mathcal O}(\log p / n)} $$ $$+ \ddot {\cal H}^{{\rm c}, {\rm c}}(0,1)  \biggl (\| x_{\bf j} \|_2^2/n +{\mathcal O} (\sqrt {\log p / n } ) \biggr ) (\bar b_{-0})_{\bf j}  
  $$
  It follows that
  $$ P_n 
(\dot g_{\bar b_0 , \bar b_{-0}, \bar d }^{{\rm b}_{-0}} )_{\bf j}  = 
P_n (\dot g_{0,0,1}^{{\rm b}_{-0}} )_{\bf j}  + \ddot {\cal H}^{{\rm c}, {\rm c}}(0,1) ( \| x_{\bf j} \|_2^2/n)  (\bar b_{-0})_{\bf j}  + {\mathcal O} ( K_{\rm x}^2 \log p / n )) .$$
Thus
$$\exp[ P_n g_{\hat b_0 , \tilde b_{-0} , \hat d } - \log \hat d ] -\exp[ P_n g_{\hat b_0 , 0 , \hat d } - \log \hat d ] $$
$$ = \exp [P_n g_{0,0,1} ] (1+ {\mathcal O} (\log p / n ))  \times $$
$$= 
\biggl ( P_n (\dot g_{0,0,1}^{{\rm b}_{-0}} )_{\bf j} ) ( \tilde b_{-0} )_{\bf j} + \ddot {\cal H}^{{\rm c}, {\rm c}}(0,1)  \| x_{\bf j} \|_2^2/n (\bar b_{-0})_{\bf j} (\tilde  b_{-0})_{\bf j}  + {\mathcal O} ( K_{\rm x}^2 \log p / n ) (\tilde b_{-0}) _{\bf j} \biggr ) $$
$$ = \underbrace{\exp [P_n g_{0,0,1} ] P_n (\dot g_{0,0,1}^{{\rm b}_{-0}} )_{\bf j} )}_{= \lambda_{\bf j}^*}  ( \tilde b_{-0} )_{\bf j}  +
\bar \kappa^2 (\| x_{\bf j} \|_2^2 / n )( \tilde b_{-0} )_{\bf j}^2 + {\mathcal O} ( K_{\rm x}^2 (\log p / n)^{3/2}) ,$$
where we used that $( \bar b_{-0} )_{\bf j} (\tilde b_{-0} )_{\bf j} \le (\tilde b_{-0} )_{\bf j}^2 $, that
$P_n (\dot g_{0,0,1}^{{\rm b}_{-0}} )_{\bf j} ) = {\mathcal O} (\sqrt {\log p / n } ) $, and that $
(\tilde b_{-0} )_{\bf j} = {\mathcal O} (\sqrt {\log p / n } ) $.

When $\lambda < \lambda_{\bf j}^*$, then $\lambda_{\bf j}^* < 0 $ implies $(\tilde b_{-0})_{\bf j} >0 $ and 
$\lambda_{\bf j}^*  > 0 $ implies   $(\tilde b_{-0})_{\bf j} < 0 $. Then
\begin{eqnarray*} 
 \lambda_{\bf j}^*  ( \tilde  b_{-0})_{\bf j}  + \lambda |(\tilde  b_{-0}) _{\bf j} |  
=  \begin{cases}   ( \lambda_{\bf j}^*+ \lambda )(\tilde  b_{-0})_{\bf j}  & (\tilde b_{-0})_{\bf j} > 0 \cr
( \lambda_{\bf j}^* - \lambda ) (\tilde  b_{-0})_{\bf j}   & (\tilde b_{-0})_{\bf j}  < 0 \cr \end{cases} 
 = - {(\lambda_{\bf j}^*  - \lambda)^2 \over   2 \bar \kappa^2 \| x_j \|_2^2 /n }.
\end{eqnarray*}
Therefore, using that $| (\tilde \beta_{0} )_{\bf j}  | = | (\tilde b_{-0})_{\bf j} |/\hat d = | (\tilde b_{-0})_{\bf j} | +
{\mathcal O} ( \log p / n )$, when $\lambda \le \lambda^* ={\mathcal O} (\sqrt {\log p / n })$,
\begin{eqnarray*}
& & \exp[ R_n ( \hat \beta_0 , \tilde \beta_{-0} , \hat \sigma) ]  + \lambda \| \tilde \beta_{-0} \|_1 
- \exp [ R_n ( \hat \beta_0, 0 , \hat \sigma) ] \\
&=&
- { \lambda_{\bf j}^*  - \lambda)^2 \over   2\bar  \kappa^2 (\| x_{\bf j} \|_2^2 /n } +
\bar \kappa^2 {(\| x_{\bf j} \|_2^2 / n )}  (\tilde b_{-0} )_{\bf j}^2  + {\mathcal O} ( K_{\rm x}^2 (\log p / n )^{3/2}   \\
& =& - { (\lambda_{\bf j}^*  - \lambda)^2 \over   4\bar  \kappa^2 \| x_j \|_2^2 /n } 
+ {\mathcal O} ( K_{\rm x}^2 (\log p / n )^{3/2}) . 
\end{eqnarray*}
If $\lambda  < \lambda^* (1- \eta) $ where $1 > \eta^2 \gg K_{\rm x}^2 \sqrt {\log p / n} $ we see that the last expression is
negative, so that $\beta_{-0}^* =0$ is not a minimizer on the set ${\cal A}$.
We get
\begin{eqnarray*}
& & \PP ( \hat \beta_{-0} \not= 0 ) \ge  \PP ( |\hat d-1| + |\hat b_0 | \ge  C \sqrt {\log p / n} )   \\
&+ &\PP ( |\hat d-1| + |\hat b_0 | \le C \sqrt {\log p / n} \wedge  \lambda^* (1- \eta) ) \ge  \lambda ) + o(1) \\
&\ge& \PP ( \lambda^* (1- \eta)\ge  \lambda ) + o(1) . 
\end{eqnarray*} 
\hfill $\sqcup \mkern -12 mu \sqcap$

\subsection{Proof of the result in Section \ref{efficiency.section}} 

{\bf Proof of Lemma \ref{efficiency.lemma}.} 
Because $\hat d $ and $\hat \beta_0$ are not penalized
$$ { \partial P_n \ell_{\beta, \sigma}  \over \partial d} \biggl \vert_{\beta= \hat \beta, d= \hat d}= 0,$$
and
$${ \partial P_n \ell_{\beta, \sigma}  \over \partial \beta_0} \biggl \vert_{\beta= \hat \beta, d= \hat d}= 0.$$
This can be rewritten as
\begin{eqnarray*}
 P_n ( \dot g_{\hat b,\hat d}^{d}) - \dot g_{\hat b,\hat d}^{{\rm b}T} (\hat \beta - \beta^*) ) &=&0 \\
P_n \hat d (g_{\hat b }^{\rm b} )_0&=& 0 .
\end{eqnarray*}

By Theorem \ref{asymptotic-continuity.theorem}, with probability tending to 1,
\begin{eqnarray*}
\biggl |  (P_n- P) ( \dot g_{\hat b,\hat d}^{d}) -\dot g_{\hat b,\hat d}^{{\rm b}T} (\hat \beta - \beta^*) - 
\dot g_{0,1}^{d})\biggr |
&=& {\mathcal O} 
(K_{\rm x} M \sqrt {\log p / n }) , \\
 \biggl | (P_n - P)  (\hat d g_{\hat b_0 , \hat b_{-0} }^{\rm b} - \dot g_{0,1}^{\rm b} ) \biggr |&= &
{\mathcal O} (
K_{\rm x} M \sqrt {\log p / n } ).
\end{eqnarray*}
But $K_{\rm x} M \sqrt {\log p / n} \stackrel{\le}{\sim}  K_{\rm x} (\lambda s^*/ \eta)  \lambda \stackrel{\le}{\sim} 
\sqrt {\lambda s^* } \lambda = o(n^{-1/2} )$, where we applied that
$\eta \gg K_{\rm x} \sqrt {\lambda s^*} $ and $\sqrt {\lambda s^{*}}  \sqrt {\log p}  \rightarrow 0$.
We also used the value of $M$ given in Theorem \ref{main2.theorem}. 
Thus with probability tending to 1, 
\begin{eqnarray*}
\biggl (P_n - P \biggr )  \begin{pmatrix}  { \partial  \ell_{\beta, \sigma}  \over \partial d} \biggl \vert_{\beta= \hat \beta, d= \hat d}\cr
 { \partial  \ell_{\beta, \sigma}  \over \partial \beta_0 } \biggl \vert_{\beta= \hat \beta , d= \hat d}\cr
\end{pmatrix}
&=& \biggl (P_n - P \biggr )  \begin{pmatrix}  { \partial  \ell_{\beta, \sigma}  \over \partial d} \biggl \vert_{\beta= \beta^*,  d= 1}\cr
 { \partial  \ell_{\beta, \sigma}  \over \partial \beta_0 } \biggl \vert_{\beta= \beta^* , d= 1}\cr
\end{pmatrix} + o(n^{-1/2} )\\
&= &P_n \begin{pmatrix}  { \partial  \ell_{\beta, \sigma}  \over \partial d} \biggl \vert_{\beta= \beta^*,  d= 1}\cr
 { \partial  \ell_{\beta, \sigma}  \over \partial \beta_0 } \biggl \vert_{\beta= \beta^* , d= 1}\cr
\end{pmatrix} + o(n^{-1/2} ).
\end{eqnarray*}
By the same arguments as used in Lemma \ref{local-expansion.lemma}, we get
$$ P  \begin{pmatrix}  { \partial  \ell_{\beta, \sigma}  \over \partial d} \biggl \vert_{\beta= \hat \beta, d= \hat d}\cr
 { \partial  \ell_{\beta, \sigma}  \over \partial \beta_0 } \biggl \vert_{\beta= \hat \beta , d= \hat d}\cr
\end{pmatrix} =
\biggl ( \ddot {\cal K} (0,1) + o(1) \biggr ) \begin{pmatrix}  \hat d-1 \cr \hat \beta_0- \beta_0^* \cr \end{pmatrix} $$
with probability tending to 1. 
 \hfill $\sqcup \mkern -12mu \sqcap$

\bibliographystyle{plainnat}
\bibliography{reference}
\end{document}